\documentclass{amsart}

\usepackage{amssymb}
\usepackage{amsmath}
\usepackage{amsthm}
\usepackage{eucal}

\usepackage{url}

\usepackage[backref,
pagebackref,%
bookmarks=true,%
colorlinks=true]%
{hyperref}

\DeclareMathOperator{\XOR}{\scriptstyle{\mathsf{XOR}}}
\DeclareMathOperator{\OR}{\scriptstyle{\mathsf {OR}}}
\DeclareMathOperator{\AND}{\scriptstyle{\mathsf {AND}}}

\DeclareMathOperator{\NOT}{\scriptstyle{\mathsf {NOT}}}
\DeclareMathOperator{\ord}{ord}
\DeclareMathOperator{\wt}{wt}

\newcommand{\B}{\mathbb B}
\newcommand{\Z}{\mathbb Z}
\newcommand{\Q}{\mathbb Q}

\renewcommand{\:}{\colon}
\renewcommand{\>}{\rightarrow}

\theoremstyle{plain}
\newtheorem{thm}{Theorem}[section]
\newtheorem{theorem}{Theorem}[section]
\newtheorem{lemma}[thm]{Lemma}
\newtheorem{proposition}[thm]{Proposition}

\newtheorem*{OpQu*}{Open question}
\theoremstyle{definition}
\newtheorem{definition}[thm]{Definition}
\theoremstyle{remark}
\newtheorem{note}[thm]{Note}
\newtheorem*{note*}{Note}
\newtheorem{example}[thm]{Example}
\newtheorem*{exmp*}{Example}
\newtheorem*{exmps*}{Examples}

\begin{document}

\title[$p$-adic Ergodicity and Pseudorandomness]{Non-Archimedean Ergodic Theory and Pseudorandom Generators}
\author{Vladimir Anashin 
}





\begin{abstract}
The paper develops techniques in order to construct computer programs, 
pseudorandom number generators (PRNG), that produce uniformly distributed sequences.
The paper exploits an approach that treats standard processor instructions (arithmetic
and
bitwise logical ones) as continuous functions on the space of 2-adic integers.
Within this approach, a PRNG is considered as a dynamical system and is
studied by means of the non-Archimedean ergodic theory.
\end{abstract}

\maketitle

\section{Introduction}
\label{Sec:Intro}

Any computer program could be viewed as a composition of basic instructions
which are the simplest instructions performed by a processor (CPU), i.e., as a composition of operators
of a proper assembler. These operators depend on a type of  CPU. 
Usually corresponding assemblers include some operators which are common for all CPUs independently of the type: these are arithmetic operators (addition, multiplication),
bitwise logical operators (e.g., $\AND$, a bitwise conjunction; $\OR$, a
bitwise disjunction, $\XOR$, a bitwise logical `exclusive or', etc.), and some others (e.g., left and right shifts). Speaking formally,
all these common operators
are defined on the set $\mathbb B^n$ of all $n$-bit words,
where
$n$ is the length of machine words the CPU operates (which is sometimes called
the CPU bitlength). However, all these common operators could be in a natural
way extended to the set $\mathbb Z_2$ of all infinite strings of zeros and
ones. The latter set $\Z_2$ could be endowed with a metric (called a 2-adic
metric) and so becomes a (non-Archimedean) metric space. Interestingly, all
these common operators are \emph{continuous functions} with respect to this
metric. So, all computer programs build from these operators could be viewed
as continuous 2-adic functions; whence, their behaviour could be studied
with the use of non-Archimedean analysis.  In this paper, we apply this
approach to construct and study pseudorandom generators.

Pseudorandom (number) generator (a PRNG for short) is a computer program that
produces a random-looking sequence of machine words, which could be also
treated as a sequence of numbers in their base-2 expansions. Pseudorandom
generators are widely used in numerous applications, especially in simulation (e.g., 
quasi Monte Carlo) and cryptography (e.g., stream ciphers). A theory (better to say, theories) of PRNG
 is an important part of computer science, see
e.g.,\cite[Chapter 3]{knuth}. We say `theorie{\bf s} of PRNG' since
the very definition of pseudorandomness assumes that the produced sequence
must pass certain class of statistical tests, so the definition of a PRNG
depends on the choice of the tests. Actually the paper could be considered
as a contribution to a non-Archimedean theory of PRNG.

As a rule, the weakest statistical property the sequence must necessarily satisfy to be
considered pseudorandom  is uniform distribution; that is, each term of the
sequence must occur with the same frequency.  For example, a well-known
linear congruential generator (LCG) 
produces the recurrence 
sequence $\{x_i\}_{i=0}^\infty$ over the set $\{0,1,\ldots,m-1\}$ 
according to the recurrence law $x_{i+1}\equiv a+bx_i\pmod m$,
for some rational integers  $a,b$. This sequence is 
uniformly distributed
if and only if it is purely periodic and the length of its shortest period is equal to the modulus $m$. The latter condition
implies that {\it each} number of $\{0,1,\ldots,m-1\}$ occurs at the period {\it exactly once} and vice versa. We refer such sequences as {\it strictly}
uniformly distributed.

In other words, the LCG produces a uniformly distributed sequence if and
only if the mapping $x\mapsto a+bx\pmod m$ of the residue ring $\Z/m\Z$ modulo
$m$ permutes residues $\{0,1,\ldots,m-1\}$ cyclically. We call the mapping
$x\mapsto a+bx$ of the ring $\Z$ of rational integers \emph{transitive modulo
$m$} in this case.  

It is not difficult to see that every composition $f$ of arithmetic and bitwise
logical operators, which defines a mapping of $\Z_2$ into $\Z_2$, induces
a well defined mapping $f\bmod 2^n$ of the residue ring $\Z/2^n\Z$ (that is, on the set $\B^n$)
into itself, for all $n=1,2,\ldots$. It turns out that the \emph{mapping $f\bmod
2^n$ is transitive for all $n$ if and only if the mapping $f$ is ergodic
(with respect to the Haar measure) on $\Z_2$}, see e.g., \cite{me-spher} for
a proof. Thus, to construct PRNGs (that produce strictly uniformly distributed
sequences over $\B^n$) out of arithmetic and bitwise logical operators we just need to construct the corresponding ergodic transformation of the space
$\Z_2$. 

This approach was already utilized in \cite{anashin2,anashin1,anashin5,anashin3,anashin4,
anashin4a,abc-v2,anashin-lect,anashin_MSU_NATO,kotomina} in order to construct numerous
\emph{non-linear} congruential generators and to study their properties.

The paper  is organized as follows:
\begin{itemize}
\item In section \ref{sec:Bas} we demonstrate that actually {a CPU works
with approximations of 2-adic integers with respect to 2-adic metric}.
\item In section \ref{sec:App} we demonstrate that both arithmetic, bitwise
logical and some other instructions of CPU could be extended to functions
that are continuous on the metric
space $\Z_2$,  
as well
as programs combined from these instructions; and that programs producing uniformly
distributed sequences could be constructed as automata  with output/state
transition functions being, accordingly, ergodic/measure preserving transformations
with respect to  a normalized Haar measure, which is a natural probabilistic measure on $\Z_2$.
\item In section \ref{sec:Tool} we develop various techniques that could be
used to construct the above mentioned ergodic/measure preserving transformations,
or to verify whether a given transformation is ergodic/measure preserving.
This section could serve mainly as a survey; however, it contains  new results
as well.
\item In section \ref{Sec:Gen} we study (with the use of the above mentioned techniques)
 two special types of fast PRNG: first one, defined by the recurrence law  $x_{i+1}\equiv
a+\sum^{m}_{j=1}a_{j}(x_i\XOR b_{j})\pmod{2^n}$, and the second one, defined
by the recurrence law $x_{i+1}\equiv a+\sum^{m }_{j=0}a_{j}\delta _{j}(x_i)$,
where $\delta_j(x)=\frac{x\AND 2^{j}}{2^{j}}$, the $j$-th binary digit in
the base-2 expansion of $x$. These generators are of special interest to
stream ciphers since they are utilized in some designs, see  \cite{abc-v2,abc_per}.
\item In section \ref{sec:Rand} we study properties of a sequence produced
by ergodic transformation of the space $\Z_2$. We demonstrate, in
particular, that this sequence satisfy D.Knuth's randomness criterion Q1,
see \cite[Section 3.5, Definition Q1]{knuth}.
\item We conclude in section \ref{sec:Conc}.
\end{itemize}

The paper is partly based on the author's preprint \cite{anashin4},
results of section \ref{Sec:Gen} were announced in author's papers \cite{anashin2,anashin1}
without proofs.

Note that most results of the paper could be re-stated 
for arbitrary
prime $p$, and \emph{not} only for $p=2$. 
%

Some  $p$-adic arguments were exploited in studies of certain special types of PRNGs, see \cite{Kl-Gor, Viv_Vlad,Wood_Smart}. However, none of these works study
PRNGs combined of basic computer instructions (both arithmetic and logical) as continuous 2-adic dynamical systems: In \cite{Kl-Gor} only an output of a feedback-with-carry
shift register  is considered as a 2-adic integer (which
actually is a rational, an irreducible fraction with odd denominator),
in \cite{Viv_Vlad,Viv_Bosk} authors study properties of pseudorandom numbers obtained
from round-off errors in calculations of 2-variate linear maps 
(actually they deal with a transformation $x\mapsto\lfloor\frac{\theta}{p^k}x\rfloor$
of the space $\mathbb Z_p$ of $p$-adic integers, where $\lfloor\cdot\rfloor$ is
an `integer part' of a $p$-adic number), in \cite{Wood_Smart} authors
study a generator with recurrence law
$x_{i+1}=\frac{x_i(x_i-1)}{2}$ on $\mathbb Z_2$, which is a
2-adic analog of a real logistic map.

It  worth noting here that there is a vast literature on PRNGs 
based on operations of finite fields and rings, see \cite{Shp_Rec_Seq} and
references therein. However, to our best knowledge none of these works use
$p$-adic techniques.

We note that the presented paper can also be considered as a contribution to
the theory of 
$p$-adic dynamical systems (especially to the $p$-adic ergodic theory). The latter theory recently attracted significant
interest due to its applications in mathematical physics, biology, genetics, cognitive sciences, etc., see e.g. \cite{khren:conf, Khren-Nils} and references
therein. However, usually  relevant works study dynamics on the whole field
$\mathbb Q_p$ of $p$-adic numbers, or even on its algebraic closure $\mathbb
C_p$, see the works just cited, as well as e.g., \cite{Ben_no_wander,Ben_hyp}.
In our paper, we study dynamical systems on $\mathbb Z_p$, which is
the ring of integers of $\mathbb Q_p$, and simultaneously  a ball of radius
1. Interestingly, our techniques developed primarily to study PRNGs was successfully applied to solve a problem (that was set up by A.~Khrennikov) on ergodicity of perturbed monomial maps on $p$-adic
spheres, see \cite{me-spher}.

%
%
\section{Basics}
\label{sec:Bas}
A contemporary processor is word-oriented.
That is, it works with words of zeroes and ones of a certain fixed length
$n$ (usually $n=8,16,32,64$). Each binary word $z\in\mathbb B^n$ of
length $n$ could be considered as a base-2
expansion of a number $z\in\{0,1,\ldots,2^n-1\}$ and vice versa.
We also can identify the set $\{0,1,\ldots,2^n-1\}$ with residues modulo
$2^n$; that is with elements of the residue ring $\mathbb Z/2^n\mathbb
Z$ modulo $2^n$. 
Actually, \emph{arithmetic} (numerical) instructions of a processor are just {\sl operations
of the residue ring} $\Z/2^n\Z$: An $n$-bit word processor performing a single instruction of addition (or multiplication)
of two $n$-bit numbers just deletes more significant digits of the sum (or
of a 
product) of these numbers thus merely reducing the result modulo $2^n$. Note
that to calculate  a sum of two integers (i.e., without reducing the result
modulo $2^n$) a `standard' processor uses not a single instruction but invokes
a program (that is a sequence of basic instructions).

Another kind of basic instructions of a processor are \emph {bitwise logical} operations:
$\XOR$, $\OR$, $\AND$, $\NOT$, which are clear from their definitions. It
worth noting only that the set $\mathbb B^n$ with respect to $\XOR$ could
be considered also as an
$n$-dimensional vector space over a field $\Z/2\Z=\mathbb B$. 

A third type of instructions could be called  \emph{machine} ones, since
they depend on the processor. But usually they include such standard instructions
as shifts (left and right) 
of an $n$-bit word.

As an example  we give formal definitions of some basic instructions (bitwise
logical and machine), the
definitions for the rest of these instructions could be obtained by an analogy.
Let 
$$z=\delta_0(z)+\delta_1(z)\cdot 2+\delta_2(z)\cdot 2^2+\delta_3(z)\cdot 2^3+\cdots$$ 
be
a base-2 expansion for $z\in\mathbb N_0=\{0,1,2,\ldots\}$ (that is, $\delta_j(z)\in\{0,1\}$).
Then, according to the respective
definitions of instructions, we have
\begin{itemize}
\item $y\XOR z=y\oplus z$ is a bitwise
addition modulo 2: $\delta_j(y\XOR z)\equiv\delta_j(y)+\delta_j(z)\pmod 2$;
\item 
$y\AND z$ is a bitwise
multiplication modulo 2: $\delta_j(y\AND z)\equiv\delta_j(y)\cdot\delta_j(z)\pmod 2$;
\item $\NOT$, 
a bitwise logical negation: $\delta_j(\NOT(z))\equiv\delta_j(z)+1\pmod 2$;
\item $\lfloor \frac{z}{2}\rfloor$, the integral part of $\frac{z}{2}$,  is a shift towards less significant
bits;
\item $2\cdot z$ is a shift towards more significant
bits;
\item $y\AND z$ is  masking of $z$ 
with the mask $y$; 
\item $z\pmod{2^k}=z\AND(2^k-1)$ is a reduction of $z$
modulo $2^k$ 
\end{itemize}
Note that in literature  $\oplus$ is used along with $\XOR$ for a bitwise `exclusive
or' operator, $\vee$ along with $\OR$, and $\wedge$ (or $\odot$) along with
$\AND$. In the rest of this paper we use only $\OR$ for bitwise logical `or',
$\AND$ for bitwise logical `and', we use $\XOR$ for `exclusive or'. 

We can make now the following important observation:
Basic instructions of a processor 
are well
defined functions on the set $\mathbb N_0$ (of  non-negative  rational integers)
valuated in $\mathbb N_0$.  

Moreover,
all mentioned basic instructions, arithmetic, bitwise logical and machine
ones, are defined on the set $\mathbb Z_2$ of all $2$-adic integers, 
which 
within the context of this paper could be thought
of as a set of all countably infinite binary sequences with terms indexed by $0,1,2,\ldots$. Sequences with only finite number of $1$s correspond to non-negative
rational integers in their base-2 expansions, sequences with only finite
number of $0$s correspond to negative rational integers, while eventually periodic 
sequences (that is, sequences that become periodic starting with a certain
place) correspond to rational numbers represented by irreducible fractions
with odd denominators: 
for instance, $3=\ldots00011$, $-3=\ldots11101$, $\frac{1}{3}=\ldots10101011$, 
$-\frac{1}{3}=\ldots1010101$. So $\delta_j(u)$ for $u\in\mathbb Z_2$ is merely
the $j$-th 
term of the corresponding sequence. 

Arithmetic operations
(addition and multiplication) with these sequences could be defined via standard
`school-textbook' algorithms  of addition and multiplication of natural numbers represented by base-2
expansions. Each term of a sequence that corresponds to the sum (respectively,
to the product) of two given sequences could be calculated 
by these algorithms with a finite number of steps. 

Thus, $\mathbb Z_2$ is
a commutative ring with respect to the so defined addition and multiplication.
It is a metric space with respect to the metric (distance) $d_2(u,v)$ defined by the following
rule:  $d_2(u,v)=\|u-v\|_2=\frac{1}{2^n}$, where $n$ is the smallest non-negative
rational integer such that
$\delta_n(u)\ne\delta_n(v)$, and $d_2(u,v)=0$ if no such $n$ exists (i.e.,
if $u=v$). For instance $d_2(3,\frac{1}{3})=\frac{1}{8}$. The function $d_2(u,0)=\|u\|_2$
is a norm of a 2-adic integer $u$, and $\ord_2 u=-\log_2\|u_2\|_2$ is a 2-adic
valuation of $u$. Note that for $u\in\mathbb N_0$ the valuation $\ord_2 u$
is merely the exponent of the highest power of 2 that divides $u$ (thus, loosely
speaking, $\ord_2 0=\infty$, so $\|0\|_2=0$).

Once the metric is defined, one defines notions of convergent sequences, limits,
continuous functions on the metric space, even derivatives
if the space is a commutative ring.  For instance, with respect to the so defined metric
on $\mathbb Z_2$ the following sequence 
tends
to $-1=\ldots 111$,
$$1,3,7,15,31,\ldots,2^n-1,\ldots\xrightarrow[d_2]{}-1,$$
bitwise logical operators (such as $\XOR,$ $\AND$, ...) define continuous
functions in two variables, the function $f(x)=x \XOR a$ is differentiable
everywhere on $\mathbb Z_2$ for every rational integer $a$: Its derivative
is $-1$ for negative $a$, and $1$ in the opposite case (see example \ref{DerLog}
for other examples of this
kind and more detailed calculations). 

Reduction
modulo $2^n$ of a $2$-adic integer $v$, i.e., setting all terms of the
corresponding sequence with indexes greater than $n-1$ to zero (that is,
taking the first $n$ digits in the representation of $v$)  is just an approximation 
of a $2$-adic integer
$v$ by a rational integer with precision $\frac{1}{2^n}$: This approximation
is an $n$-digit positive rational integer $v \AND (2^n-1)$; the latter will
be denoted also as $v\bmod{2^n}$.

Actually {\it a processor works
with approximations of 2-adic integers with respect to 2-adic metric}:
When
one tries to load a number whose base-2 expansion contains more than $n$
significant bits into a registry of an $n$-processor, the processor just
writes only $n$ low order bits of the number in a registry { thus reducing
the number modulo} $2^n$. Thus, 
precision of the approximation is defined by the bitlength of the processor.

All these considerations (after proper modifications) remain true for arbitrary
prime $p$, and \emph{not} only for $p=2$, thus leading to the notion of a $p$-adic integer
and to  $p$-adic analysis. For formal introduction to $p$-adic
analysis, exact notions and results see any relevant book, e.g. \cite{Kobl,Mah}.

\section{Approach}
\label{sec:App}

Arithmetic and bitwise logical operations are not independent:
Some of them could be expressed via the others. For instance, for all
$u,v\in\mathbb Z_2$
\begin{equation}
\label{eq:id}
\begin{split}
&\NOT u=u \XOR (-1);\\
&u+\NOT u=-1;\\
&u \XOR v = u+v-2(u \AND v);\\ 
&u \OR v = u+v-(u \AND v);\\
&u \OR v=(u \XOR v)+(u \AND v).
\end{split}
\end{equation}
Proofs of   
identities \eqref{eq:id} are just an exercise: For example, if 
$\alpha ,\beta\in\{0,1\}$ then $\alpha\XOR \beta=\alpha+\beta -2\alpha\beta$
and $\alpha\OR  \beta=\alpha+\beta -\alpha\beta$. Hence:
\begin{multline*}
u \XOR v =\sum_{i=0}^\infty 2^i(\delta_i (u)\XOR \delta_i (v))=\\
\sum_{i=0}^\infty 2^i(\delta_i (u)+\delta_i (v)-2\delta_i (u)\delta_i (v))=\\
\sum_{i=0}^\infty 2^i(\delta_i (u))+\sum_{i=0}^\infty 2^i(\delta_i (v))-
2\cdot\sum_{i=0}^\infty 2^i(\delta_i (u)\delta_i (v))=\\
u+v-2(u \AND v).
\end{multline*}
Proofs of the remaining identities can be made by analogy
and thus are omitted. A shift towards more
significant digits, as well as masking 
could
be derived from the above operations: An $m$-step shift of $u$ is $2^m u$;
masking of $u$ is $u \AND M$, where $M$ is an integer which base-2 expansion
is a mask (i.e., a string of $0$s and $1$s). 

A common feature  the above mentioned arithmetic, bitwise logical and machine
operations all share is that they are, with the only exception of shifts towards less significant
bits, {\it compatible}, that is,  $\omega(u,v)\equiv\omega(u_1,v_1)\pmod{2^r}$
whenever both congruences $u\equiv u_1\pmod{2^r}$  and $v\equiv v_1\pmod{2^r}$ hold
simultaneously (here $\omega$ stands for any of these operations, arithmetic, bitwise
logical,  or machine). The notion of a compatible mapping could be naturally generalized
to mappings $(\mathbb Z/2^l\mathbb Z)^{t}\rightarrow(\mathbb Z/2^l\mathbb
Z)^{s}$ and 
$\mathbb Z_2^{t}\rightarrow\mathbb Z_2^{s}$ of Cartesian products.

We note that considerations we made above,  after  proper modifications hold for arbitrary prime $p$, and not only for $p=2$. The case of odd prime
$p$ is important to produce pseudorandom sequences on $N$ symbols, $N>2$. 
PRNGs that produce pseudorandom numbers
in the range $\{0,1,2,\ldots, N-1\}$ are often used in practice,
and we are going to discuss them also. However, the case $p=2$ will be sometimes
exceptional in our considerations (this often happens in $p$-adic analysis),
so from time to time we have to switch to the case $p=2$ and then revert
back to the general case. 

The compatibility property, being originally stated in algebraic terms, could be expressed in terms of $p$-adic analysis
as well, for arbitrary prime $p$, and not only for $p=2$. Namely this is not difficult to verify that the \emph{function $F\colon\mathbb Z_p^{t}\rightarrow\mathbb Z_p^{s}$ is compatible if and only if it satisfies
Lipschitz condition with coefficient 1 with respect to $p$-adic distance};
e.g., for $s=t=1$ the function $F$ is compatible if and only if
$$
\|F(u)-F(v)\|_p\le \|u-v\|_p
$$
for all $u,v\in\mathbb Z_p$.

Obviously, a composition of compatible mappings
is a compatible mapping. We list now some important examples of compatible
operators $(\mathbb Z_p)^{t}\rightarrow(\mathbb Z_p)^{s}$, $p$ prime. 
Here are some of them that originate from arithmetic operations:

\begin{equation}
\begin{split}
& {\text {\rm multiplication,}}\ \cdot:\ (u,v)\mapsto uv;\\ 
& {\text {\rm addition,}}\ +:\ (u,v)\mapsto u+v; 
\\
& {\text {\rm subtraction,}}\ -:\ (u,v)\mapsto u-v;
\\
& {\text {\rm exponentiation,}}\ \uparrow_p:\ (u,v)\mapsto u\uparrow_p v=(1+pu)^v;
\\ 
& {\text {\rm raising to negative powers}},\ u\uparrow_p(-n)=(1+pu)^{-n}; 
\\ 
& {\text {\rm division,}}\ /_p: u/_pv=u\cdot (v\uparrow_p(-1))=\frac{u}{1+pv}. 
\label{eq:opAr}
\end{split}
\end{equation}

The other part originates from digitwise logical operations of $p$-valued logic:
\begin{equation}
\label{eq:opLog} 
\begin{split}
& {\text {\rm digitwise multiplication}}\ u\odot_p v:\\ 
& \delta_j(u\odot_p
v)\equiv \delta_j (u)\delta_j (v)\pmod p;\\ 
& {\text {\rm digitwise addition}}\ 
u\oplus_p v:\\ 
& \delta_j(u\oplus_p
v)\equiv \delta_j (u)+\delta_j (v)\pmod p;\\ 
& {\text {\rm digitwise subtraction}}\
u\ominus_p v:\\ 
& \delta_j(u\ominus_p
v)\equiv \delta_j (u)-\delta_j (v)\pmod p.
\end{split}
\end{equation}
Here 
$\delta_j(z)$ $( j=0,1,2,\ldots)$
stands for the $j$-th 
digit of $z$ in its base-$p$ expansion. For $p=2$ equations \eqref{eq:opLog}
define $\AND$ and $\XOR$.  
%

In case $p=2$ compatible mappings could be characterized in terms of Boolean
functions. Namely, each transformation $T\colon\mathbb Z/2^n\Z\rightarrow\mathbb Z/2^n\Z$ of the residue ring $Z/2^n\Z$ modulo $2^n$
could be
considered as an ensemble of $n$ Boolean functions $\tau_i^T(\chi_0,\ldots,\chi_{n-1})$,
$i=0,1,2,\ldots,n-1$,
in $n$ Boolean variables $\chi_0,\ldots,\chi_{n-1}$ by assuming $\chi_i=\delta_i(u)$,
$\tau_i^T(\chi_0,\ldots,\chi_{n-1})=\delta_i(T(u))$
for $u$ running from $0$ to $2^n-1$. The following easy proposition
holds.
\begin{proposition}
\label{Bool} 
{\rm \cite
{anashin2}}
A mapping $T\colon\mathbb Z/2^n\Z\rightarrow\mathbb Z/2^n\Z$ 
{\rm (}accordingly, a mapping $T\colon\mathbb Z_2\rightarrow\mathbb Z_2${\rm
)}
is compatible
if and only if each Boolean function $\tau_i^T(\chi_0,\chi_{1},\ldots)=\delta_i(T(u))$,
$i=0,1,2,\ldots$,
does not depend on the variables $\chi_{j}=\delta_j(u)$ for $j>i$.
\end{proposition}
\begin{note*}
We use the term `compatible' instead of the
term `conservative' of \cite{anashin2}, since the latter term in numerous papers
on algebraic systems has attained another meaning, see  \cite[p. 45]{LN}.
Note that in the theory of Boolean functions mappings satisfying conditions of the proposition are also known
as {\it triangular} mappings, and  as \emph{T-functions} in cryptography. 

The proposition after proper restatement (in
terms of functions of $p$-valued logic) also holds for 
odd prime $p$. For multivariate mappings proposition \ref{Bool} holds also:
a mapping $T=(t_1,\ldots,t_s)\colon\mathbb Z_2^{r}\rightarrow\mathbb 
Z_2^{s}$ is compatible
if and only if each Boolean function $\tau_i^{t_j}(\chi_{1,0},\chi_{1,1},\ldots,
\chi_{r,0},\chi_{r,1},\ldots)=\delta_i(t_k(u,\ldots,u_r))$ ($i=0,1,2,\ldots$,
$k=0,1,\ldots,s$) does not depend on variables $\chi_{\ell,j}=\delta_j(u_{\ell})$ 
for $j>i$ ($\ell=1,2,\ldots,r$).
\end{note*}

Now, given a compatible mapping $T\colon\mathbb Z_2\rightarrow\mathbb Z_2$, one
can define an
induced mapping 
$T\bmod2^n\colon\mathbb Z/2^n\Z\rightarrow\mathbb Z/2^n\Z$
assuming $(T\bmod 2^n)(z) =T(z)\bmod 2^n=(T(z))\AND(2^n-1)$ 
for $z=0,1,2,\ldots,2^n-1$. The induced mapping is obviously a
compatible mapping of the ring $\mathbb Z/2^n\Z$ into itself. For odd prime
$p$, as well as for multivariate case 
$T\colon\mathbb Z_p^{s}\rightarrow\mathbb Z_p^{t}$
an induced mapping $T\bmod p^n$ could be defined by analogy.
\begin{definition}
\label{def:erg}
We call a compatible mapping $T\colon\mathbb Z_p\rightarrow\mathbb Z_p$
{\it bijective modulo $p^n$} if and only if the induced mapping $T\bmod p^n$ is a permutation
on $\mathbb Z/p^n\Z$; we call $T$ {\it transitive modulo $p^n$}, if and only if $T\bmod p^n$
is a  permutation with a single cycle. 
We call a compatible mapping
$T\colon\mathbb Z_p^{s}\rightarrow\mathbb Z_p^{t}$
{\it balanced modulo $p^n$} if and only if the induced mapping $T\bmod p^n$ maps
$(\mathbb Z/p^n\Z)^{s}$ onto $(\mathbb Z/p^n\Z)^{t}$, and each element of 
$(\mathbb Z/p^n\Z)^{t}$ has the same number of preimages in $(\mathbb Z/p^n\Z)^{s}$.
\end{definition}

Often a pseudorandom generator 
could be constructed as a  finite automaton
${\mathfrak A}=\langle N,M,f,F,u_0\rangle $ with a finite state set $N$, state
transition function
$f:N\rightarrow N$, finite output alphabet $M$, output function  $F:N\rightarrow M$
and an initial state (seed) $u_0\in N$. The following sequence $\mathcal T=\{u_j=f^j(u_0)\}_{j=0}^\infty$
is called a sequence of states:
\begin{equation*}
f^{j}(u_0)=\underbrace{f(\ldots f(}_{j
\;\text{times}}u_0)\ldots)\ \ (j=1,2,\ldots); f^{0}(u_0)=u_0.
\end{equation*}
Thus, the generator produces the output sequence $\mathcal S$ over
the set $M$ out of the sequence of states:
\begin{equation*}
\mathcal S=F(u_0), F(f(u_0)), F(f^{2}(u_0)),\ldots, F(f^{j}(u_0)),\ldots
\end{equation*}

Mappings that are transitive modulo $p^n$, as well as mappings that are  balanced modulo $p^n$  could
be used as building blocks of pseudorandom generators to provide both large
period
length and uniform distribution of output sequences. Namely, the following obvious
proposition holds.
\begin{proposition}
\label{prop:Auto}
If the state transition function $f$ of the automaton $\mathfrak A$ is
transitive on the state set $N$, i.e., if $f$ is a permutation with a single cycle
of length $|N|$, if, further, $|N|$ is a multiple of $|M|$, and if the output function 
$F:N\rightarrow M$ is balanced
{\rm (}i.e.,  $|F^{-1}(s)|=|F^{-1}(t)|$ for all $s,t\in M${\rm )}, then the output sequence
$\mathfrak S$ of the automaton $\mathfrak A$ is purely periodic with period length  $|N|$ 
{\rm (i.e., maximum possible)}, and each element of  
$M$ occurs at the period the same number of times, $\frac{|N|}{|M|}$ exactly. {\rm
That
is, the
output sequence $\mathcal S$ is strictly uniformly distributed.} 
\end{proposition}
Note that in case $N=\mathbb B^{kn}$ and $M=\mathbb B^{ln}$ one can use a transitive modulo
$2^{kn}$ compatible state transition function $f\colon\Z/2^{kn}\Z\rightarrow\Z/2^{kn}\Z$ and a balanced 
modulo $2^n$ output function $F\colon(\Z/2^{n}\Z)^k\rightarrow (\Z/2^{n}\Z)^l$ to produce
a strictly uniformly distributed sequence.
  
Now we describe connections between generators of strictly uniformly distributed
sequences and $p$-adic ergodic theory.
Recall that a {\it dynamical system on a measurable space} $\mathbb
S$ is a triple $(\mathbb S;\mu; f)$, where $\mathbb S$ is a set endowed with a measure
$\mu$, and $f\colon \mathbb S\> \mathbb S$ is a {\it measurable function}; that is, an
$f$-preimage of any measurable subset is a measurable subset. These basic
definitions from dynamical system theory, as well as the following ones,
could be found at \cite{KN}; see also \cite{Kat_Has} as a comprehensive monograph on various aspects of
dynamical systems theory.

A {\it trajectory} of a dynamical system is a sequence
$$x_0, x_1=f(x_0),\ldots, x_i=f(x_{i-1})=f^i(x_0),\ldots$$
of points of the space $\mathbb S$, $x_0$ is called an {\it initial} point of the
trajectory. If $F\colon \mathbb S\> \mathbb T$ is a measurable mapping to some other measurable space
$\mathbb T$ with a measure $\nu$ (that is, if an $F$-preimage of any $\nu$-measurable subset
of $\mathbb T$ is a $\mu$-measurable subset of $X$), the sequence $F(x_0), F(x_1), F(x_2),\ldots$ is called an {\it observable}. Note that the trajectory formally
looks like the sequence of states of a pseudorandom generator, whereas the
observable resembles the output sequence.

A mapping $F\colon\mathbb S\rightarrow\mathbb
Y$ of a measurable space  $\mathbb S$ into a measurable space $\mathbb Y$ endowed with probabilistic measure  $\mu$ and $\nu$, respectively,
is said to be 
{{\it measure preserving}} (or, sometimes, {\it equiprobable})
whenever $\mu(F^{-1}(S))=\nu(S)$ for each measurable subset   $S\subset\mathbb Y$. In case $\mathbb S=\mathbb Y$ and $\mu=\nu$,
a measure preserving mapping  $F$ is said to be 
{{\it ergodic}} whenever for each measurable subset 
$S$ such that $F^{-1}(S)=S$ holds either  $\mu(S)=1$ or $\mu(S)=0$.

Recall
that to define a measure $\mu$  on some set $\mathbb S$ we should assign non-negative real numbers
to some subsets that are called elementary. All other {\it measurable} subsets
are compositions of these elementary subsets with respect to countable unions,
intersections, and complements. 

Elementary subsets in $\Z_p$ are balls $B_{p^{-k}}(a)=a+p^k\Z_p$ of radii
$p^{-k}$ (in other
words,
co-sets with respect to ideal generated by $p^k$). To each ball we assign a number $\mu_p(B_{p^{-k}}(a))=\frac{1}{p^k}$. This way we define
a probabilistic measure on the space $\Z_p$, 
$\mu_p(\Z_p)=1$. The measure $\mu_p$
is called a (normalized) {\it Haar measure} on $\Z_p$. The normalized Haar
measure on $\mathbb Z_p^n$ could be defined by analogy.

Note that the sequence $\{s_i\}_{i=0}^\infty$ of $p$-adic
integers is uniformly distributed (with respect to the normalized Haar measure
$\mu_p$   on $\mathbb
Z_p$) 
if and only if it is uniformly distributed modulo $p^k$ for all $k=1,2,\ldots$;
That is, for every $a\in\mathbb Z/p^k\Z$ relative numbers of occurrences 
of $a$ in the initial segment of length $\ell$ in the sequence 
$\{s_i\bmod p^k\}$ of residues
modulo $p^k$ 
are asymptotically equal,
i.e., 
$\lim_{\ell\to\infty}\frac{A(a,\ell)}{\ell}=\frac{1}{p^k}$, where 
$A(a,\ell)=|\{s_i\equiv a\pmod{p^k}\colon i<\ell\}|$, see \cite{KN} for
details. 
Thus, strictly uniformly distributed sequences are uniformly distributed
in the common sense of theory of distributions of sequences.
Moreover, the following theorem (which was announced in \cite{anashin3} and
proved in \cite{me-spher}) holds.
\begin{theorem} 
\label{thm:erg-tran}
For $m=n=1$, a compatible mapping  $F\colon\mathbb Z_p^n\rightarrow\mathbb
Z_p^m$ 
{preserves the normalized Haar measure} $\mu_p$ on $\mathbb Z_p$ 
\textup{(}resp., 
{is ergodic} with respect to $\mu_p$
\textup{)}
if and only if it is 
{
bijective 
\textup{(}resp., 
transitive
\textup{)} modulo $p^k$} for all $k=1,2,3,\ldots$  

For $n\ge m$,  the mapping $F$ {preserves measure} $\mu_p$ if and only if it induces a {balanced} mapping of $(\mathbb Z/p^k\Z)^n$ onto $(\mathbb Z/p^k\Z)^m$, for all 
$k=1,2,3,\ldots$.

\end{theorem}
This theorem in combination with proposition \ref{prop:Auto} implies in particular that whenever one chooses a compatible and ergodic mapping $f\colon\Z_2\rightarrow\Z_2$  as a state transition function
of the automaton $\mathfrak A$, and a compatible and measure-preserving mapping
$F\colon(\Z/2^{n}\Z)^k\rightarrow (\Z/2^{n}\Z)^l$ as an output function of
$\mathfrak A$, both the sequence of states and output sequence of the
automaton are uniformly distributed with respect to the Haar measure. This
implies that reduction of these sequences  modulo $2^n$ results in strictly
uniformly distributed sequences of binary words. 
Note also that reduction modulo $2^n$ a computer performs automatically.

Thus, theorem \ref{thm:erg-tran} gives us a way to construct generators
of uniformly distributed sequences out of standard computer instructions. 
%
%
%
Now the problem is how to describe these measure preserving (in particular,
ergodic) mappings in the class of all compatible mappings. We start to develop some
theory to answer the following questions: What compositions of  {basic  instructions}
are measure preserving? are ergodic? Given a composition of basic instructions, is it measure preserving? is it ergodic?

%
%

\section{Tools}
\label{sec:Tool}
In this section we introduce various techniques in order to construct
measure preserving and/or ergodic mappings, as well as to verify whether
a given mapping is measure preserving or, respectively, ergodic. We are
mainly focused on the class of compatible mappings.

Main results of Subsection \ref{ssec:Interp} are Theorem \ref{ergBin} and
Theorem \ref{ergAnGen}. With the use of these
one can verify whether a given function 
is 
measure-preserving, or ergodic. Theorem \ref{ergBin} gives a general method
yet demands a function must be represented via interpolation series. Theorem \ref{ergAnGen} gives an easier method for a narrower class of functions,
which is, however, rather wide: e.g., it contains polynomials and rational functions.

The main result of Subsection \ref{subsec:comb} is Theorem \ref{Delta}, which
gives a general method how to construct a measure-preserving or ergodic fucntion
out of arbitrary compatible function.

Theorem \ref{ergBool} is the central point of Subsection \ref{subsec:Bool}.
Being more of theoretical value, it has as a consequence a useful Proposition
\ref{compBool}, which gives an easy method to construct new vast classes
of ergodic functions out of given ergodic function.

Subsection \ref{subsec:unider} deals with differentiation. In particular,
this subsection
introduces Calculus for functions build from basic computer operators. 
The main result of this
subsection is Theorem \ref{ergDer} which gives  conditions for a uniformly differentiable
function to be ergodic.

\subsection{Interpolation series}
\label{ssec:Interp}

The general characterization of compatible ergodic functions is given by
the following theorem.
\begin{theorem}
\label{ergBin}
{\rm{\cite{anashin2,anashin1}}} 
A function $f\colon{\mathbb Z}_{2}\rightarrow {\mathbb Z}_{2}$ 
is compatible if and only if it can be represented as 
\begin{equation*}
f(x)=c_0+\sum^{\infty }_{i=1}c_{i}\,2^{\lfloor \log_2 i \rfloor}\binom{x}{i}
\qquad (x\in\mathbb Z_2);
\end{equation*}
The function $f$ 
is compatible and measure preserving if and only if it can be represented as 
\begin{equation*}
f(x)=c_0+x+\sum^{\infty }_{i=1}c_{i}\,2^{\lfloor \log_2 i \rfloor +1}\binom{x}{i}
\qquad (x\in\mathbb Z_2);
\end{equation*}
The function 
$f$
is compatible and ergodic if and only if it can be represented as 
\begin{equation*}
f(x)=1+x+\sum^{\infty}_{i=1}c_{i}2^{\lfloor \log_{2}(i+1)\rfloor+1}\binom{x}{i}
\qquad (x\in\mathbb Z_2),
\end{equation*}
where $c_0,c_1, c_2 \ldots \in {\mathbb Z}_2$.
\end{theorem}
Here, as usual, 
\begin{equation*}
\binom{x}{i}=
\begin{cases}\dfrac{x(x-1)\cdots  (x-i+1)}{i!}, 
& \text {for $i=1,2,\ldots$};\cr 1, & \text{for $i=0$},
\end{cases}
\end{equation*}
and $\lfloor\alpha\rfloor$ is the integral part of $\alpha$, i.e.,
the largest rational integer not exceeding $\alpha$.
%
\begin{note*} 
For odd prime $p$ an analog of the statement of theorem \ref{ergBin}
provides
only sufficient conditions for ergodicity (resp., measure preservation)
of $f$: namely, if $(c,p)=1$,
i.e., if $c$ is a unit (=invertible element) of $\mathbb Z_p$, then the
function
$f(x)=c+x+\sum^{\infty}_{i=1}c_{i}p^{\lfloor \log_{p}(i+1)\rfloor+1}\binom{x}{i}$ 
defines a compatible and ergodic mapping of $\mathbb Z_p$ onto itself,
and 
the
function
$f(x)= c_0+c\cdot x+\sum^{\infty}_{i=1}c_{i}p^{\lfloor \log_{p}i\rfloor+1}\binom{x}{i}$ 
defines a compatible and measure preserving mapping of $\mathbb Z_p$ onto itself
(see \cite{anashin3}).
\end{note*}
Thus, in view of theorem \ref{ergBin} one can choose a state transition
function to be a polynomial with rational (not necessarily integer)
coefficients setting $c_i=0$ for all but finite number of $i$.
Note that to determine whether a given polynomial $f$ with rational (and not
necessarily integer) coefficients is integer valued (that is, maps $\mathbb
Z_p$ into itself), compatible and ergodic, it is sufficient to determine
whether it
induces a permutation with a single cycle of $O(\deg f)$ integral points. To be more exact, the following
proposition holds.
\begin{proposition}
\label{prop:Qpol} 
{\rm\cite{anashin3}}
A polynomial $f(x)\in {\Q}_{p}[x]$ over the field of $p$-adic numbers $\Q_p$
is integer valued, compatible, and ergodic
{\rm (}resp., measure preserving{\rm)} if and only if 
$$z\mapsto f(z)\bmod p^{\lfloor
\log_p (\deg f)\rfloor +3},$$ 
where $z$ 
runs through $0,1,\ldots,p^{\lfloor
\log_p (\deg f)\rfloor +3}-1$,  is a compatible and transitive 
{\rm (}resp., bijective{\rm)} mapping
of the residue ring $\Z/p^{\lfloor
\log_p (\deg f)\rfloor +3}\Z$ onto itself. 
\end{proposition}

Although
this is not very essential for further considerations, we note, however, that
the series in the statement of theorem \ref{ergBin} and of the note thereafter are
uniformly convergent with respect to $p$-adic distance. Thus
the mapping $f\colon\mathbb Z_p\rightarrow\mathbb Z_p$ is well defined
and continuous with respect to $p$-adic distance, 
see 
\cite[Chapter 9]{Mah}.

Theorem \ref{ergBin} can be applied in design of \emph{exponential} (the ones based
on exponentiation) generators  of uniformly distributed sequences.  
\begin{example} 
\label{expGen}
For any odd $a=1+2m$ the function $f(x)=ax+a^x$  
is transitive modulo $2^n$, for all $n=1,2,\ldots$

Indeed, in view of theorem \ref{ergBin} the function $f$ defines a compatible and ergodic
transformation of $\mathbb Z_2$ 
since $f(x)=(1+2m)x+(1+2m)^x=x+2mx+\sum_{i=0}^\infty m^i 2^i\binom{x}{i}=
1+x+4m\binom{x}{1}+
\sum_{i=2}^\infty m^i 2^i\binom{x}{i}$ and $i\ge\lfloor\log_2(i+1)\rfloor+1$
for all $i=2,3,4,\ldots$. 

This generator could be of practical value since it uses not more than
$n+1$ multiplications modulo $2^n$ of $n$-bit numbers; of course, one should
use calls to the look-up table 
$a^{2^j}\bmod{2^n}$, $j=1,2,3,\ldots,n-1$. The latter table must be precomputed,
 corresponding calculations involve $n-1$ multiplications modulo $2^n$. 
 
\end{example}
\begin{note*}
A similar argument shows that for every prime
$p$ and  every $a\equiv 1\pmod
p$ the function $f(x)=ax+a^x$ defines a compatible and ergodic mapping
of $\mathbb Z_p$ onto itself.
\end{note*}

For polynomials with  (rational or $p$-adic) integer coefficients 
theorem \ref{ergBin} may be restated in the following form.
\begin{proposition}
\label{ergPol}
{\rm \cite{anashin2, anashin1}}
Represent a polynomial $f(x)\in\mathbb Z_2[x]$  in a basis of descending 
factorial powers
$$
x^{\underline 0}=1,\ x^{\underline 1}=x,
\ldots,x^{\underline i}=x(x-1)\cdots(x-i+1),\ldots,$$
that is, let
$$f(x)=\sum^{d}_{i=0}c_i\cdot x^{\underline i}$$
for $c_0,c_1,\dots,c_d\in\mathbb Z_2$. Then the polynomial $f$ induces
an ergodic {\rm (and, obviously, a compatible)} mapping of $\mathbb Z_2$ onto
itself if and only if its coefficients $c_0,c_1,c_2, c_3$ satisfy the following congruences:   
$$
\begin{aligned}
c_0\equiv & 1\pmod 2,\ c_1\equiv 1\pmod 4,\\ 
c_2\equiv & 0\pmod 2,\ c_3\equiv 0\pmod 4.
\end{aligned}
$$
The polynomial $f$ induces a measure preserving mapping if and only if
$$c_1\equiv 1\ (\bmod\, 2),\quad c_2\equiv 0\ (\bmod\, 2),\quad c_3\equiv 0\ (\bmod\, 2).$$
\end{proposition}
Thus, to provide ergodicity of the polynomial  $f$  
it is necessary and sufficient to fix $6$ bits only, while the other bits of 
coefficients of $f$ may be arbitrary.  This guarantees transitivity
of the state transition function $z\mapsto f(z)\bmod 2^n$ for each $n$, and hence,
uniform distribution of the sequence of states.

Proposition \ref{ergPol} implies that the polynomial $f(x)\in\mathbb
Z[x]$ is ergodic (resp., measure preserving) if and only if it is transitive modulo 8
(resp., if and only if it is bijective modulo 4). A corresponding assertion
holds in a general case, for arbitrary prime $p$.
\begin{theorem}
\label{ergPolGen}
{\rm \cite{larin}} A polynomial $f(x)\in\mathbb Z_p[x]$
induces an ergodic transformation of $\mathbb Z_p$ if and only if it is transitive
modulo $p^2$ for $p\ne 2,3$, or modulo $p^3$, for $p=2,3$. The polynomial
$f(x)\in\mathbb Z_p[x]$ induces a measure preserving transformation of 
$\mathbb Z_p$ if and only if it is bijective
modulo $p^2$.
\end{theorem}

\begin{example} 
The mapping $x\mapsto f(x)\equiv x+2x^2\pmod{2^{32}}$ (which is used in
a cipher RC6, see \cite{RC6}) is bijective, since it is bijective modulo 4: 
$f(0)\equiv 0\pmod4$, $f(1)\equiv 3\pmod4$, $f(2)\equiv 2\pmod4$, 
$f(3)\equiv 1\pmod4$. Thus, the mapping 
$x\mapsto f(x)\equiv x+2x^2\pmod{2^{n}}$ is bijective for all $n=1,2,\ldots$. 
\end{example}
Hence, with the use of theorem \ref{ergPolGen} it is possible to
construct transitive modulo $q>1$ mappings for arbitrary natural $q$: One
just takes $f(z)=(1+z+\hat qg(z))\bmod q$, where $g(x)\in\mathbb Z[x]$ is
an arbitrary polynomial, and $\hat q$ is a product of $p^{s_p}$ for all
prime factors $p$ of $q$, where $s_2=s_3=3$, and $s_p=2$ for $p\ne 2,3$. For example, a polynomial $f(x)=201+201x+200x^{17}$ is transitive modulo $10^n$
for arbitrary $n$.  

In these considerations,
the polynomial $g(x)$ may be chosen, roughly speaking, `more or less at random',
yet the output sequence will be uniformly
distributed for any choice of $g(x)$. This assertion can be generalized
also:
\begin{proposition}
\label{ergAn} {\rm \cite{anashin3}} Let
$p$ be a prime, and let
$g(x)$ be an arbitrary composition of arithmetic operations \textup{(see
\eqref{eq:opAr} of section \ref{sec:App})}.
Then the mapping $z\mapsto 1+z+p^2g(z)$\ $(z\in\mathbb Z_p)$ is ergodic. 
\end{proposition}

In fact, both propositions \ref{ergPol}, \ref{ergAn} and theorem \ref{ergPolGen} 
are 
special cases of the following general theorem.
\begin{theorem}
\label{ergAnGen}
{\rm \cite{anashin3}}
Let $\mathcal B_p$ be a class of all functions defined by series of
the form $f(x)=\sum^{\infty}_{i=0}c_i\cdot x^{\underline i}$, where 
$c_0,c_1,\dots$ are $p$-adic integers, and
$x^{\underline i}$, $i=0,1,2,\ldots$,  
are descending factorial powers {\rm(see proposition \ref{ergPol})}.
Then
the function $f\in \mathcal B_p$ preserves measure if and only if it is bijective
modulo $p^2$; $f$ is ergodic if and only if it is transitive modulo $p^2$
{\rm(}for $p\ne 2,3${\rm)}, or modulo $p^3$ {\rm(}for $p\in\{2,3\}${\rm)}.

\end{theorem}
\begin{note*} As it was shown in \cite{anashin3}, the class $\mathcal B_p$
contains all polynomial functions over $\mathbb Z_p$,
as well as analytic (e.g., rational, entire) functions that are convergent everywhere
on $\mathbb Z_p$. Actually, every mapping that is  a composition 
of arithmetic operators \eqref{eq:opAr}  belong to $\mathcal B_p$; thus, every
such mapping modulo $p^n$ could be induced by a polynomial with rational
integer coefficients (see the end of Section 4 in \cite{anashin3}). For instance,
the mapping $x\mapsto (3x+3^x) \bmod 2^n$ (which is transitive modulo $2^n$,
see example \ref{expGen}) could be induced by the polynomial $1+x+4\binom{x}{1}+
\sum_{i=2}^{n-1} 2^i\binom{x}{i}=1+5x+\sum_{i=2}^{n-1} \frac{2^i}{i!} \cdot
x^{\underline i}$ --- just note that $c_i=\frac{2^i}{i!}$ are $2$-adic integers
since the exponent of maximal power of $2$ that is a factor of $i!$ 
is exactly $i-\wt_2i$,
where $\wt_2 i$ is a number of $1$s in the base-2 expansion of $i$
(see e.g. \cite[Chapter 1, Section 2, Exercise 12]{Kobl}); thus 
$\|c_i\|_2=2^{-\wt_2 i}\le 1$, i.e. $c_i\in \mathbb Z_2$ and so $c_i\bmod
{2^n}\in \mathbb Z$.  
\end{note*}

Theorem \ref{ergAnGen} implies that, for instance, the state transition
function $f(z)=(1+z+\zeta(q)^2(1+\zeta(q)u(z))^{v(z)})\bmod q$ is transitive
modulo $q$ for each natural $q>1$ and arbitrary polynomials $u(x),v(x)\in\mathbb
Z[x]$,  where $\zeta(q)$ is a product of all prime factors of $q$. So 
one can choose as a state transition function not only polynomial functions, but
also rational functions, as well as analytic ones. 
For instance, certain \emph{inversive generators} (that exploit multiplicative inverses of residues modulo $2^n$) could
be considered. 
\begin{example}
\label{Invers}
The function $f(x)=-\frac{1}{2x+1}-x$ 
is transitive modulo $2^n$, for all $n=1,2,3,\ldots$.

Indeed, the function $f(x)=(-1+2x-4x^2+8x^3-\cdots)-x=-1+x-4x^2+8(\cdots)$ is analytic
and is defined everywhere on $\mathbb Z_2$; thus $f\in\mathcal B_p$. Now the
conclusion follows from theorem \ref{ergAnGen} since by direct calculations
it could be easily verified that the function $f(x)\equiv -1+x-4x^2\pmod
8$ is transitive modulo 8. Note that 
the mapping $x\mapsto
f(x)\bmod 2^n$ could be induced by the polynomial $-1+x-4x^2+8x^3+\cdots+(-1)^n
2^{n-1}x^{n-1}$.
\end{example}
\subsection{Combinations of operators}
\label{subsec:comb}
A transformation  of the residue ring  $\Z/q\Z$ induced by a polynomial with rational integer coefficients is the only type of mapping that could be constructed
as a composition of arithmetic operations, $+$ and $\cdot$.
The class of all transitive modulo $q$ mappings induced by polynomials with
rational integer coefficients is rather wide: For instance, for $q=2^n$
it contains $2^{O(n^2)}$ mappings (for exact value 
see  \cite[Proposition 16]{larin}). 
However, this class could
be widened significantly (up to a class of order $2^{2^n-n-1}$ in case
$q=2^n$) by  including bitwise logical operators 
into the composition. Actually, every compatible mapping could be constructed
this way.

\begin{proposition}
\label{erg-comp} 
Let 
$g$
be a compatible
mapping of $\mathbb Z_2$ onto itself. Then for each $n=1,2,\ldots$
the mapping $\bar g=g\bmod 2^n$ could be represented as a finite
composition of arithmetic and bitwise logical operators \textup{(actually,
as a composition of $+$, $\XOR$, $\AND$ and shifts towards higher order bits,
i.e., multiplications by powers of 2) }. 
\end{proposition}
\begin{proof}
In view of proposition \ref{Bool}, one could
represent $\bar g$ as 
\begin{multline*}
\bar g(x)=\gamma_0(\chi_0)+2\gamma_1(\chi_0,\chi_1)+\cdots \\
+2^{n-1}\gamma_{n-1}(\chi_0,\ldots,\chi_{n-1}),
\end{multline*}
where $\gamma_i=\delta_i(\bar g)$, $\chi_i=\delta_i(x)$, $i=0,1,\ldots,n-1$.
Since each $\gamma_i(\chi_0,\ldots,\chi_i)$ is a Boolean function in Boolean
variables $\chi_0,\ldots,\chi_i$, it
could be expressed via finite number of $\XOR$s and $\AND$s of these variables
$\chi_0,\ldots,\chi_i$. Yet each variable $\chi_j$ could be expressed as
$\chi_j=\delta_j(x)=2^{-j}(x\AND(2^j))$; thus 
\begin{multline*}
2^i\gamma_i(\chi_0,\ldots,\chi_i)=
\gamma_i(2^{i}(x\AND(1)),2^{i-1}(x\AND(2)),\ldots\\ 2(x\AND(2^{i-1})),x\AND(2^{i})),
\end{multline*}
and the conclusion follows. 
\end{proof}
  
It turns out that there is an \emph{easy way to construct a measure preserving or ergodic
mapping out of an arbitrary compatible mapping:} 
\begin{theorem}
\label{Delta} \textup{\cite{anashin3}} Let $\Delta$ be a difference operator, i.e., $\Delta g(x)=g(x+1)-g(x)$
by definition. Let, further, $p$ be a prime, let $c$ be coprime with
$p$, $\gcd(c,p)=1$, and let $g\colon\mathbb Z_p\rightarrow
\mathbb Z_p$ be a compatible mapping. Then the mapping $z\mapsto c+z+p\Delta
g(z)\ (z\in\mathbb Z_p)$ is ergodic, and the mapping $z\mapsto d+cz+pg(z)$, 
preserves measure for arbitrary $d$.

Moreover, if $p=2$, then the converse also holds: Each compatible and ergodic
\textup {(}respectively, each compatible
and measure preserving\textup {)}
mapping $z\mapsto f(z)\ (z\in\mathbb Z_2)$ can be represented as
$$f(x)=1+x+2\Delta g(x)$$  \textup {(}respectively as
$f(x)=d+x+2g(x)$\textup {)} for suitable $d\in\mathbb Z_2$ and compatible 
$g\colon\mathbb
Z_2\rightarrow \mathbb Z_2$.
\end{theorem}
\begin{note*} The case $p=2$ is the only case where the converse of the first
assertion of the proposition \ref{Delta} holds. 
\end{note*}
\begin{example}
\label{KlSh-2} 
Proposition \ref{Delta} immediately implies
Theorem 2 of \cite{KlSh}: For any composition $f$ of primitive functions, 
the mapping $x\mapsto x +2f(x )\pmod {2^n}$ is invertible --- just note
that
a composition of primitive
functions is compatible (see \cite{KlSh} for the definition of primitive
functions).\qed
\end{example}
Theorem \ref{Delta} could be an important tool in design of pseudorandom
generators, since it provides high flexibility during design. 
In fact, one may use nearly
arbitrary composition of arithmetic and bitwise logical operators to  produce a
strictly uniformly distributed sequence:
Both for 
$g(x)=x\XOR(2x+1)$ and for 
$$g(x)=\Biggl(1+2\frac{x\AND x^{2}+x^3\OR x^4}{3 + 4(5+6x^5)^{x^6\XOR x^7}}\Biggr)^{7+\frac{8x^8}{9+10x^9}}$$   
(note, both these functions $g$ are compatible!) the sequence $\{x_i\}$ defined by the recurrence relation
$x_{i+1}=(1+x_i+2(g(x_i+1)-g(x_i)))\bmod {2^n}$ is strictly uniformly distributed
in $\mathbb Z/2^n\Z$,  for all $n=1,2,3\ldots$. 
Actually, a designer could vary the function $g$ in a very
wide scope without worsening prescribed values of 
some important statistical characteristics of output sequence. 
As a matter of fact,  choosing proper  arithmetic and bitwise logical operators 
the designer is restricted
only by desirable performance since any compatible ergodic mapping could
be produced  this way.

\subsection{Boolean representation}
\label{subsec:Bool}
In case $p=2$ the two preceding  subsections give  two (equivalent) complete
descriptions of the class of all
compatible ergodic mappings, namely, theorem \ref{ergBin} and theorem
\ref{Delta}. They enable one to express {\it any} compatible and transitive
modulo $2^n$ state transition function either as a polynomial of special kind
over a field
$\mathbb Q$ of rational numbers, or as a special composition of arithmetic
and bitwise logical operations. 
Both these representations are suitable for programming, since they involve
only standard machine instructions. However, we
need one more representation, in a Boolean form (see proposition \ref{Bool}). Although this
representation is not very convenient for programming,
it outlines 
some new methods for construction
of ergodic transformations, see proposition \ref{compBool} below. Also, this
representation could be of use while proving the ergodicity of some
simple
mappings, see e.g. example \ref{KlSh-3} below.
The following theorem is just a restatement of a known (at least 30 years
old) result from the theory 
of Boolean functions, the so-called bijectivity/transitivity criterion  for triangle
Boolean mappings.
However, the latter is mathematical folklore, and thus it is somewhat
difficult to
attribute it, yet a reader can find a proof in, e.g.,
\cite[Lemma 4.8]{anashin2}.
\begin{theorem}
\label{ergBool} 
A mapping $T\colon\mathbb Z_2\rightarrow\mathbb Z_2$ is
compatible and measure preserving if and only if for each $i=0,1,\ldots$ the algebraic normal form, ANF, of the Boolean function 
$\tau^T_i=\delta_i(T)$
in Boolean variables $\chi_0,\ldots,\chi_{i}$ can be represented as 
 
$$\tau^T_i(\chi_0,\ldots,\chi_i)=\chi_i+\varphi^T_i(\chi_0,\ldots,\chi_{i-1}),$$ 
where $\varphi^T_i$ is an ANF
of a Boolean function in Boolean variables $\chi_0,\ldots,\chi_{i-1}$. The mapping $T$ is compatible  and ergodic if and only if,
in addition to already stated conditions, the following conditions hold:
$\varphi^T_0=1$, and each Boolean function
$\varphi^T_i$ $(i>0)$ is of odd weight. 
\end{theorem}
\begin{note*} 
Recall that the {\it algebraic normal form} (ANF for short) of the Boolean function $\psi(\chi_0,\ldots,\chi_j)$
is the representation of this function via $\oplus$ (addition modulo 2, that
is,
logical `exclusive or') and $\odot$ (multiplication modulo 2, that is, logical `and', or
conjunction). In other words, the ANF of the Boolean function $\psi$ is its representation
in the form 
\begin{multline*}
\psi(\chi_0,\ldots,\chi_j)= \beta\oplus\beta_0\odot\chi_0\oplus\beta_1\odot\chi_1\oplus\ldots\\
\oplus\beta_{0,1}\odot\chi_0\odot\chi_1\oplus\ldots,
\end{multline*}
where $\beta,\beta_0,\ldots\in\{0,1\}$. The ANF is sometimes called a {\it
Boolean polynomial}. In the sequel in the ANF we write $+$ instead of $\oplus$ and $\cdot$
instead of $\odot$ when this does not lead to misunderstanding.


Recall that  {\it weight} of the Boolean function $\psi$ in $(j+1)$
variables is the number of $(j+1)$-bit words that { satisfy} $\psi$;
that is, weight of a Boolean function is cardinality of a truth set of the Boolean function.

Note that \emph{weight of the Boolean function $\varphi (\chi_0,\ldots,\chi_{i-1)}$
in Boolean  variables $\chi_0,\ldots,\chi_{i-1}$ is odd
if and only if degree $\deg \varphi$ of the Boolean function
$\varphi$   is exactly
$i$, that is, if and only if the ANF of  $\varphi$ contains a monomial
$\chi_0\cdots\chi_{i-1}$.}

\end{note*}


\begin{example}
\label{KlSh-3} 
With the use of  theorem \ref{ergBool} it is possible to give a short proof of the
main result of \cite{KlSh}, namely, of Theorem 3 there:
{\it The mapping $f (x)=x +(x^2\OR C )$ over $n$-bit words is invertible
if and only if the least significant bit of $C$ is 1. For $n\ge 3$ it is a permutation
with a single cycle if and only if both the least significant bit and the third least
significant bit of $C$ are $1$.}

{\it Proof of theorem 3 of} \cite{KlSh}. 
Recall that for $x\in\mathbb Z_2$ and $i=0,1,2,\ldots$ 
we denote $\chi_i=\delta_i(x)\in\{0,1\}$; also we denote $c_i=\delta_i(C)$. 
We will calculate ANF of the Boolean function $\delta_i(x+(x^2\OR C))$  
in variables $\chi_0,\chi_1,\ldots$. We start with the following easy claims:
\begin{itemize}
\item $\delta_0(x^2)=\chi_0$,\ $\delta_1(x^2)=0$,\ $\delta_2(x^2)=\chi_0\chi_1+\chi_1$,
\item $\delta_n(x^2)=\chi_{n-1}\chi_0+\psi_{n}(\chi_0,\ldots,\chi_{n-2})$ for all
$n\ge 3$, where $\psi_{n}$ is a Boolean function in $n-1$
Boolean
variables $\chi_0,\ldots,\chi_{n-2}$.
\end{itemize}

The first of these claims could be easily verified by direct calculations. To prove
the second one represent $x=\bar x_{n-1}+2^{n-1}s_{n-1}$ for $\bar x_{n-1}=x\bmod
2^{n-1}$ and
calculate $x^2=(\bar x_{n-1}+2^{n-1}s_{n-1})^2=
\bar x_{n-1}^2+2^{n}s_{n-1}\bar x_{n-1}+2^{2n-2}s_{n-1}^2=\bar x_{n-1}^2+2^n\chi_{n-1}\chi_0
\pmod{2^{n+1}}$ for $n\ge 3$ and note that $\bar x_{n-1}^2$ depends only on 
$\chi_0,\ldots,\chi_{n-2}$.

This gives
\begin{enumerate}
\item $\delta_0(x^2\OR C)=\chi_0+c_0+\chi_0c_0$
\item $\delta_1(x^2\OR C)=c_1$
\item $\delta_2(x^2\OR C)=\chi_0\chi_1+\chi_1+c_2+c_2\chi_1+c_2\chi_0\chi_1$
\item $\delta_n(x^2\OR C)=\chi_{n-1}\chi_0+\psi_{n}+c_n+c_n\chi_{n-1}\chi_0+c_n\psi_{n}$
for $n\ge 3$
\end{enumerate} 
From here it follows  
that if $n\ge 3$, then $\delta_n(x^2\OR C)=
\lambda_n(\chi_0,\ldots,\chi_{n-1})$, 
and $\deg \lambda_n\le
n-1$, since $\psi_{n}$ depends only on 
$\chi_0,\ldots,\chi_{n-2}$.

Now we successively calculate $\gamma_n=\delta_n(x+(x^2\OR C))$ for $n=0,1,2,\ldots$.
We have $\delta_0(x+(x^2\OR C))=c_0+\chi_0c_0$ so necessarily $c_0=1$
since otherwise $f$ is not bijective modulo 2. Proceeding further with
$c_0=1$ we obtain $\delta_1(x+(x^2\OR C))=c_1+\chi_0+\chi_1$, since
$\chi_1$ is a carry. Then $\delta_2(x+(x^2\OR C))=(c_1\chi_0+c_1\chi_1+\chi_0\chi_1)+
(\chi_0\chi_1+\chi_1+c_2+c_2\chi_1+c_2\chi_0\chi_1)+\chi_2=
c_1\chi_0+c_1\chi_1+\chi_1+c_2+c_2\chi_1+c_2\chi_0\chi_1+\chi_2$,
here $c_1\chi_0+c_1\chi_1+\chi_0\chi_1$ is a carry. From here in view of \ref{ergBool}
we immediately deduce that $c_2=1$ since otherwise $f$ is not transitive
modulo 8. 
Now for $n\ge 3$ 
one has $\gamma_n=\alpha_{n}+\lambda_n
+\chi_n$, where $\alpha_n$ is a carry, and $\alpha_{n+1}=\alpha_n\lambda_n
+\alpha_n\chi_n+\lambda_n\chi_n$. But if $c_2=1$ then $\deg\alpha_3=\deg
(\mu\nu+\chi_2\mu+\chi_2\nu)=3$, where $\mu=c_1\chi_0+c_1\chi_1+\chi_0\chi_1$, 
$\nu=(\chi_0\chi_1+\chi_1+c_2+c_2\chi_1+c_2\chi_0\chi_1)=
0$. This implies inductively in view of (iv) above that 
$\deg\alpha_{n+1}=n+1$ and that $\gamma_{n+1}=\chi_{n+1}+\xi_{n+1}(\chi_0,\ldots,\chi_{n})$,
$\deg\xi_{n+1}=n+1$. So conditions of \ref{ergBool} are satisfied, thus
finishing the proof of theorem 3 of \cite{KlSh}.\qed
\end{example}

There are some other applications of Theorem \ref{ergBool}.
\begin{proposition}
\label{compBool}
Let 
$F\colon\mathbb Z_2^{n+1}\rightarrow\mathbb Z_2$ be a compatible mapping
such that for all $z_1,\ldots,z_n\in\mathbb
Z_2$ the mapping $F(x,z_1,\ldots,z_n)\colon\mathbb Z_2\rightarrow \mathbb
Z_2$ is  measure preserving. Then $F(f(x),2g_1(x),\ldots,2g_n(x))$ preserves
measure for all compatible $g_1,\ldots,g_n\colon\mathbb Z_2\rightarrow \mathbb
Z_2$ and all compatible and measure
preserving $f\colon\mathbb Z_2\rightarrow \mathbb
Z_2$. Moreover, if 
$f$ is ergodic then $f(x+4g(x))$, $f(x\XOR (4g(x)))$, $f(x)+4g(x)$, and
$f(x)\XOR  (4g(x))$ are ergodic for any compatible $g\colon\mathbb Z_2\rightarrow \mathbb
Z_2$
\end{proposition}
\begin{proof} Since the function $F$ is compatible, $\delta_i(F(u_0,u_1,\ldots,u_n)$
does not depend on $\delta_j(u_k)=\chi_{j,k}$ for $j>i$ (see proposition
\ref{Bool}
and a note thereafter). Consider ANF of the Boolean function $\delta_i(F(u_0,u_1,\ldots,u_n))$:
\begin{multline*}
\delta_i(F(u_0,u_1,\ldots,u_n))=\\
\chi_{0,i}\Psi_i(u_0,u_1,\ldots,u_n)+
\Phi_i(u_0,u_1,\ldots,u_n),
\end{multline*}
where Boolean functions 
$\Psi_i(u_0,u_1,\ldots,u_n)$ and $\Phi_i(u_0,u_1,\ldots,u_n)$ do not depend on 
$\chi_{0,i}$; 
that is, they depend only on 
$$\chi_{0,0},\ldots,\chi_{0,i-1},
\chi_{1,0},\ldots,\chi_{1,i},\ldots,\chi_{n,0},\ldots,\chi_{n,i}.$$
In view of theorem \ref{ergBool}, $\Psi_i=1$ since $F(x,z_1,\ldots,z_n)$ 
preserves measure for all $z_1,\ldots,z_n\in\mathbb Z_2$. Moreover, then
$\Phi_i(f(x),2g_1(x),\ldots,2g_n(x))$ does not depend on $\chi_i=\delta_i(x)$
since $\delta_j(2g(x))$ does not depend on $\chi_i$ for all $j=1,2,\ldots,n$.
So in view of theorem \ref{ergBool},
$\delta_i(f(x))=\chi_i+\xi_i(f(x))$, where $\xi_i(f(x))$ does not depend
on $\chi_i$ since $f$ preserves measure. 
Finally,
\begin{multline*} 
\delta_i(F(f(x),2g_1(x),\ldots,2g_n(x)))=\\
\delta_i(f(x))+
\Phi_i(f(x),2g_1(x),\ldots,2g_n(x))=\\
\chi_i + \xi_i(f(x))+
\Phi_i(f(x),2g_1(x),\ldots,2g_n(x))=\chi_i+\Xi_i,
\end{multline*} 
where 
the Boolean function
$\Xi_i$ depends only on 
$\chi_0,\ldots,\chi_{i-1}$. 
This proves the 
first assertion of proposition \ref{compBool} in view of theorem \ref{ergBool}.

We prove the second assertion along similar lines.  For $z\in\mathbb Z_2$
and $i=0,1,2,\ldots$
let $\zeta_i=\delta_i(z)$. Thus one can represent $\delta_i(z\XOR 4g(z))$ and
$\delta_i(z+ 4g(z))$ via 
ANFs  in Boolean variables $\zeta_0,\zeta_1,\ldots,\zeta_i$. Note that
$\delta_i(z\XOR 4g(z))=\zeta_i+\lambda_i(z)$, where $\lambda_i(z)=0$
for $i=0,1$ and $\deg\lambda_i(z)\le i-1$ for $i>1$, since for $i>1$ the Boolean
function $\lambda_i(z)$ depends 
only on $\zeta_0,\ldots,
\zeta_{i-2}$. 

Further, we claim that 
$\delta_i(z+ 4g(z))=\delta_i(z)+\mu_i(z)$, where 
$\mu_i(z)=\mu_i^g(z)$ is 0 for $i=0,1$ and $\deg\mu_i(z)\le i-1$ for $i>1$. 
Indeed, $\mu_i(z)=\lambda_i(z)+\alpha_i(z)$, where the Boolean function
$\alpha_i(z)$ is a carry. Yet $\alpha_i(z)=0$ for $i=0,1,2$, 
and 
$\alpha_i(z)=\zeta_{i-1}\lambda_{i-1}(z)+\zeta_{i-1}\alpha_{i-1}(z)+
\lambda_{i-1}(z)\alpha_{i-1}(z)$ for $i\ge 3$, and $\alpha_i(z)$ depends
only on 
$\zeta_0,\ldots,\zeta_{i-1}$ since $\alpha_i(z)$
is a carry. However, $\deg\alpha_3(z)=2$ and if $\deg\alpha_{i-1}(z)\le
i-2$ then $\deg\delta_{i-1}(z)\alpha_{i-1}(z)\le i-1$, 
$\deg\lambda_{i-1}(z)\alpha_{i-1}(z)\le i-1$, and 
$\deg\zeta_{i-1}\lambda_{i-1}(z)\le i-1$ since $\alpha_{i-1}(z)$ depends
only on 
$\zeta_0,\ldots,\zeta_{i-2}$ and
$\lambda_{i-1}(z)$ depends 
only on $\zeta_0,\ldots,
\zeta_{i-3}$. Thus $\deg\alpha_i(z)\le i-1$ and hence $\deg\mu_i(z)\le
i-1$.

Now, since $f(x)$ is
ergodic, $\delta_{i}(f(x))=\chi_i+\xi_i(x)$, where the Boolean function
$\xi_i$ depends only on 
$\chi_0,\ldots,\chi_{i-1}$ and,
additionally, $\xi_0=1$, and $\deg\xi_i=i$ for $i>0$ (see theorem \ref{ergBool});
i.e. $\xi_i(x)=\chi_0\chi_1\cdots\chi_{i-1}+\vartheta_i(x)$,
where $\deg\vartheta_i(x)\le i-1$ for $i>0$.
Hence, for $\ast\in\{+,\XOR\}$ one has
$\delta_{i}(f(x\ast 4g(x)))=\delta_i(x\ast 4g(x))+
\delta_0(x\ast 4g(x))\delta_1(x\ast 4g(x))\cdots\delta_{i-1}(x\ast 4g(x))+
\vartheta_i(x\ast 4g(x))$; thus $\delta_{i}(f(x\ast 4g(x)))=\chi_i+
\chi_0\cdots\chi_{i-1}+ \beta_i^\ast(x)$, where $\deg\beta_i^\ast(x)\le
i-1$ for $i>0$, and $\delta_0(f(x\ast 4g(x))=\delta_0(x\ast 4g(x))+1=\chi_0+1$.
Finally, $f(x\ast 4g(x))$ for $\ast\in\{+,\XOR\}$ is ergodic in view of
theorem
\ref{ergBool}.

In a similar manner it could be demonstrated that $f(x)\ast 4g(x)$ is ergodic
for $\ast\in\{+,\XOR\}$: $\delta_i(f(x)\ast 4g(x))=\delta_i(f(x))$ 
for $i=0,1$ and thus satisfy the conditions of theorem \ref{ergBool}. For $i>1$
one has $\delta_i(f(x)\XOR 4g(x))=\chi_i+\xi_i(x)+\delta_{i-2}(g(x))$;
but $\delta_{i-2}(g(x))$ does not depend on $\chi_{i-1},\chi_{i}$.
Thus the Boolean function $\xi_i(x)+\delta_{i-2}(g(x))$ in variables
$\chi_0,\ldots,\chi_{i-1}$ is of odd weight, since $\xi_i(x)$
is of odd weight, thus proving that $f(x)\XOR 4g(x)$ is ergodic.

Now  represent $g(x)=g(f^{-1}(f(x)))=h(f(x))$, where $f^{-1}$ is the
inverse mapping for $f$. Clearly, $f^{-1}(x)$ is well defined since 
the mapping $f\colon\mathbb Z_2\rightarrow\mathbb Z_2$ is bijective; 
moreover $f^{-1}(x)$
is compatible and ergodic. Finally  
$\delta_i(f(x)+ 4g(x))=\delta_i(f(x))+\mu_i^\prime(f(x))$, 
where the ANF of the  Boolean function
$\mu_i^\prime(x)=\mu_i^{h}(x)$ in Boolean variables 
$\chi_0,\ldots,\chi_{i-1}$ does not contain a monomial 
$\chi_0\cdots\chi_{i-1}$ (see the claim above). This implies that the ANF
of the Boolean function
$\mu_i^\prime(f(x))$ in Boolean variables $\chi_0,\ldots,\chi_{i-1}$
does not contain a monomial $\chi_0\cdots\chi_{i-1}$ either,
since $\delta_j(f(x))=\chi_j+\xi_j(x)$ and $\xi_j(x)$ depend only
on $\chi_0,\ldots,\chi_{j-1}$ for $j=2,3,\ldots$. Hence,
$\delta_i(f(x)+ 4g(x))=\chi_i+\xi_i(x)+\mu_i^\prime(f(x))$ and the
Boolean function $\xi_i(x)+\mu_i^\prime(f(x))$ in Boolean variables
$\chi_0,\ldots,\chi_{i-1}$ is of odd weight. This finishes the
proof in view of theorem \ref{ergBool}. 
\end{proof}
\begin{example}
\label{XOR}
With the use of \ref{compBool} it is possible to construct very fast generators
$x_{i+1}=f(x_i)\bmod 2^n$ that are transitive modulo $2^n$. 
For instance, take
$$f(x)=(\ldots((((x+c_0)\XOR d_0)
\cdots +c_m)\XOR
d_m,$$
where $c_0\equiv 1\pmod 2$, and the rest of $c_i,d_i$ are 0 modulo 4.
In a  general situation these functions $f$ (for arbitrary $c_i, d_i$) 
were studied in \cite{kotomina}, where it was proved that  $f$
is ergodic if and only if it is transitive modulo 4.
\end{example} 

\subsection{Uniform differentiability}
\label{subsec:unider}
In previous subsections we consider some methods that could be used to verify
whether a given transformation $f$ of the space $\Z_2$ is measure preserving
or ergodic. One way is to represent $f$ by interpolation series and apply
theorem \ref{ergBin}, the second way is to represent $f$ in a special form
described by theorem \ref{Delta}, the third way is to use Boolean representation
and theorem \ref{ergBool}. These methods are universal meaning they could
be applied to any compatible function $f$. However, they work only in a univariate
case.

In this subsection we present another method that works for multivariate
functions also, but is not universal any more; the method could be applied
only to uniformly differentiable mappings and some mappings that are close
to these. The class of these mappings is rather wide, though.

Now we recall a generalized version of 
the main notion of
Calculus, a derivative {\it modulo} $p^k$, which was originally introduced
in \cite{anashin2,anashin1,anashin3}.  By  the definition, for points $ \mathbf a=(a_{1},\ldots  ,a_{n})$ and
 $\mathbf b=(b_{1},\ldots  ,b_{n})$ of
${\mathbb Z}^{n}_{p}$ 
the congruence $\mathbf  a\equiv  \mathbf b\pmod{p^{s}}$ means that
$\|a_{i}-b_{i}\| _{p}\le p^{-s}$ (or, the same, that $a_{i}=b_{i}+c_{i}p^{s}$ 
for suitable $c_{i}\in {\mathbb Z}_{p}$, $i=1,2,\ldots  ,s$); that is $\|\mathbf
a-\mathbf b\|_p\le p^{-s}$.

\begin{definition}
\label{def:Der}
A function $$F=(f_{1},\ldots  ,f_{m})\colon{\mathbb Z}^{n}_{p}\rightarrow {\mathbb Z}^{m}_{p}$$
is said to be {\it differentiable modulo $p^k$} at the point 
$ \mathbf u=(u_{1},\ldots  ,u_{n})\in {\mathbb Z}^{n}_{p}$
if there exists a positive
integer rational
$N$ and an $n\times m$ matrix $F^{\prime}_{k}(\mathbf u)$ over ${\mathbb Q}_{p}$
\textup{(called {\it the Jacobi matrix modulo} $p^{k}$ of the function $F$ at the
point
$\mathbf u$)} such that for every positive rational integer 
$K\ge N$ and every $ \mathbf h=(h_{1},\ldots  ,h_{n})\in {\mathbb Z}^{n}_{p}$ the congruence 
\begin{equation}
\label{Der} 
F( \mathbf u+\mathbf h)\equiv F(\mathbf u)+ 
\mathbf hF^{\prime}_{k}(\mathbf u)\pmod{p^{k+K}}
\end{equation}
holds whenever $\|\mathbf h\| _{p}\le     p^{-K}$.
 In case $m=1$ the
Jacobi matrix modulo $p^k$ is called a {\it differential modulo $p^k$}. In
case $m=n$ a determinant of the Jacobi matrix modulo $p^k$ is called a {\it Jacobian
modulo $p^k$}. Entries of the Jacobi matrix modulo $p^k$
are called {\it partial derivatives modulo} $p^k$ of the function $F$ at
the point $\mathbf u$.
\end{definition} 
A partial derivative (respectively, a differential) modulo $p^k$ is
sometimes  denoted as 
$\frac{\partial_k f_i (\mathbf u)}{\partial_k x_j}$ (respectively, as
$d_{k}F(\mathbf u)=\sum^n_{i=1} \frac {\partial_k F(\mathbf u)}{\partial_k x_i}d_{k}x_{i}$).

Since the notion of function that is differentiable modulo $p^k$ is of high
importance for the theory that follows, we discuss this notion in detail. 
Compared to differentiability, the {differentiability modulo $p^k$} is a weaker restriction. Speaking loosely, in a univariate case ($m=n=1$), definition \ref{def:Der} just yields
that   
$$
\frac{F( \mathbf u+\mathbf h)- F(\mathbf u)}{\mathbf h}\approx F^{\prime}_{k}(\mathbf u)$$
Note that whenever $\approx$ (`approximately') stands for an `{\it arbitrarily high} precision' one obtains a common definition of differentiability; however,
if  $\approx$ stands for a `precision that is {\it not worse than} $p^{-k}$',
one obtains the differentiability modulo $p^k$.

We note that the notion of  a derivative modulo $p^k$ have no direct analog in the classical
Calculus: A derivative with a precision up to the $k$-th digit after the
point, being often used in common speech, is meaningless from
the rigorous point of view since there is no distinguished  base in real
analysis. However, this notion is meaningful in $p$-adic analysis since there
is a distinguished base; namely, base-$p$. 

In $p$-adic analysis, it is obvious that whenever a function is differentiable (and its derivative
is a $p$-adic integer), it is differentiable
modulo $p^k$ for all $k=1,2,\ldots$, and in this case the derivative modulo
$p^k$ is just a {\sl reduction} of a derivative modulo $p^k$ (note that according
to definition \ref{def:Der} partial derivatives modulo $p^k$ are determined
up to a summand that is 0 modulo $p^k$).

 In cases when all partial derivatives
modulo $p^k$ at all points of  
$\mathbb Z_p^{n}$ are
$p$-adic integers, we say that the function 
$F$ has {\it integer valued derivative modulo} $p^k$; 
in these cases we can associate to each partial derivative modulo $p^k$
a unique element of the ring $\mathbb Z/p^k\Z$; 
a Jacobi matrix modulo $p^k$ 
at each point $\mathbf u\in \mathbb Z_p^{n}$ 
thus can be considered as a matrix over a ring $\mathbb Z/p^k\Z$. It turns
out that this is exactly the case for a compatible function $F$. Namely, the following
proposition holds.
\begin{proposition}
\label{intDer}
\textup{\cite{anashin2,anashin1}}
Let a compatible function 
$F=(f_{1},\ldots  ,f_{m})\colon{\mathbb Z}^{n}_{p}\rightarrow {\mathbb Z}^{m}_{p}$ be uniformly
differentiable modulo $p^k$ at the point $\mathbf u\in {\mathbb Z}^{n}_{p}$.
Then $\big\|\frac{\partial_k f_i (\mathbf u)}{\partial_k x_j}\big\|_p\le 1$, i.e.,
$F$ has integer valued derivatives modulo $p^k$.
\end{proposition}
For functions with integer valued derivatives modulo $p^k$  
the `rules of differentiation
modulo $p^k$' have the same (up to congruence modulo $p^k$ instead of equality)
form as for usual differentiation.
For instance, if both functions 
$G\colon{\mathbb Z}^{s}_{p}\rightarrow {\mathbb Z}^{n}_{p}$ and
$F\colon{\mathbb Z}^{n}_{p}\rightarrow {\mathbb Z}^{m}_{p}$ 
are differentiable modulo 
$p^{k}$ at the points, respectively, $\mathbf v=(v_{1},\ldots  ,v_{s})$
and $\mathbf u=G(\mathbf v)$, and their partial derivatives modulo $p^{k}$ at
these points are $p$-adic integers, then a composition 
$F\circ G\colon{\mathbb Z}^{s}_{p}\rightarrow {\mathbb Z}^{m}_{p}$ 
of these functions is uniformly differentiable modulo $p^{k}$ at the point
$\mathbf v$, all its partial derivatives 
modulo $p^{k}$ at this point are $p$-adic integers, and 
$(F\circ G)^\prime_k (\mathbf v)\equiv G^\prime_k (\mathbf v) F^\prime_k (\mathbf u)\pmod
{p^k}$.

\begin{definition}
\label{def:uniDer}
A function $F\colon{\mathbb Z}^{n}_{p}\rightarrow {\mathbb Z}^{m}_{p}$ is said
to be 
{\it uniformly differentiable modulo $p^k$ on $\mathbb Z_p^{n}$} if and only if there
exists $K\in\mathbb N$ such that congruence \eqref{Der} holds simultaneously for all 
$\mathbf u \in \mathbb Z_p^{n}$ as soon as
$\| h_{i}\| _{p}\le     p^{-K}$, $(i=1,2,\ldots  ,n)$. The
least of these
$K$
is denoted $N_k(F)$. 
\end{definition}
Recall that   \emph {all  partial derivatives
 modulo $p^k$ of a uniformly differentiable modulo $p^k$ function $F$
 are periodic functions with period  
$p^{N_k(F)}$}, see \cite[Proposition  2.12]{anashin2}. 
Thus,  \emph{each partial derivative modulo
$p^k$ could be considered as a function defined on (and valuated in) the
residue ring $\mathbb Z/p^{N_k(F)}\Z$}. 
Moreover, if a continuation $\tilde F$ of the function
$F=(f_{1},\ldots , f_{m})\colon{\mathbb N}^{n}_{0}\rightarrow {\mathbb N}^{m}_{0}$  
to the space $\mathbb Z_p^{n}$ is a uniformly differentiable modulo $p^k$ function on 
$\mathbb Z_p^{n}$, then one can simultaneously continue the function $F$ together with  all its
(partial) derivatives modulo $p^k$ to the whole space $\mathbb Z_p^{n}$. Consequently, we may study if necessary (partial) 
derivatives modulo $p^k$
of the function $\tilde F$ instead of those of $F$ and vise versa.
For example, a partial derivative $\frac{\partial_k f_i (\mathbf u)}{\partial_k x_j}$
modulo $p^k$ vanishes modulo $p^k$ at no point of  $\mathbb Z_p^{n}$
(that is,
$\frac{\partial_k f_i (\mathbf u)}{\partial_k x_j}\not\equiv 0\pmod{p^k}$
for all $u\in \mathbb Z_p^{(n)}$, or, the same
$\big\|\frac{\partial_k f_i (\mathbf u)}{\partial_k x_j}\big\|_p> p^{-k}$
everywhere on $\mathbb Z_p^{n}$) if and only if 
$\frac{\partial_k f_i (\mathbf u)}{\partial_k x_j}\not\equiv 0\pmod{p^k}$
for all $u\in\{0,1,\ldots,p^{N_k(F)}-1\}$.

In case $p=2$, differentiation modulo $p^k$ could naturally be implemented as a computer
program since this differentiation just implies (for a univariate $F$) estimation of the fraction $\frac{F( \mathbf u+\mathbf h)-F(\mathbf h)}{\mathbf h}$ with a $k$-bit precision, i.e., evaluation
of the first $n$ low order bits of the base-2 expansion of the corresponding
number. 
To calculate a derivative of, for instance, a state transition function, which
is a composition of basic instructions of CPU (that is, of `elementary' functions, see proposition \ref{erg-comp})
one needs to know derivatives of these  `elementary' functions, 
such as arithmetic and bitwise logical operations. Here we briefly introduce a $p$-adic analog of  
a `table of derivatives' of a classical Calculus.
\begin{example} Derivatives of bitwise logical operations.
\label{DerLog}
\begin{enumerate}
\item {\it a function $f(x)=x\AND c$ is uniformly differentiable on $\mathbb
Z_2$ for any $c\in
\mathbb Z$; $f^\prime(x)=0$ for $c\ge 0$, and $f^\prime(x)=1$ for $c<0$,} since
$f(x+2^ns)=f(x)$, and 
$f(x+2^ns)=f(x)+2^ns$ for $n\ge l(|c|)$, where $l(|c|)$ is the bit length
of absolute value of $c$
(mind that for $c\ge 0$ the $2$-adic representation
of $-c$ starts with $2^{l(c)}-c$ in less significant bits followed by $\ldots
11$:
$-1=\ldots 111$, $-3=\ldots 11101$, etc.).
\item {\it a function $f(x)=x\XOR c$ is uniformly differentiable on $\mathbb
Z_2$ for any $c\in
\mathbb Z$; $f^\prime(x)=1$ for $c\ge 0$, and $f^\prime(x)=-1$ for $c<0$.} This
immediately
follows from (i) since $u\XOR v=u+v-2(x\AND v)$ (see \eqref{eq:id} in section
\ref{sec:App}); thus
$(x\XOR c)^\prime=x^\prime+c^\prime-2(x\AND c)^\prime=1+2\cdot(0,\ \text{for}\
c\ge 0;\ \text{or}\
-1,\ \text{for}\ c<0)$.
\item in the same manner it could be shown that {\it functions $(x\bmod
2^n)$, $\NOT(x)$
and $(x\OR c)$ for $c\in \mathbb Z$ are uniformly differentiable on $\mathbb
Z_2$, and $(x\bmod 2^n)^\prime=0$, $(\NOT x)^\prime=-1$, 
$(x\OR c)^\prime=1$ for $c\ge 0$, 
$(x\OR c)^\prime=0$ for $c< 0$.}
\item {\it a function $f(x,y)=x\XOR y$ is not uniformly differentiable on 
$\mathbb Z_2^{2}$ (as a bivariate function),
yet it is uniformly differentiable modulo $2$ on $\mathbb Z_2^{2}$};
from (ii) it follows that its partial derivatives modulo 2 are 1 everywhere
on $\mathbb Z_2^{2}$.
\end{enumerate}
\end{example} 

Here is how it works altogether.
\begin{exmps*}
A function $f(x)=x+(x^2\OR 5)$ is uniformly differentiable
on $\mathbb Z_2$, 
and $$f^\prime (x)=1+2x\cdot
(x\OR 5)^\prime=1+2x.$$

A function $$F(x,y)=(f(x,y),g(x,y))=
(x \XOR 2(x \AND y ),(y +3 x^3 )\XOR x )$$ 
is uniformly differentiable modulo $2$ as a bivariate
function, and $N_1(F)=1$; namely
\begin{multline*}
F(x+2^nt,y+2^ms)\equiv\\
 F(x,y)+(2^nt,2^ms)\cdot
\begin{pmatrix}
1&x+1\\
0&1
\end{pmatrix}
\pmod{2^{k+1}}
\end{multline*}
for all $m,n\ge 1$ (here $k=\min\{m,n\}$). The matrix
$\begin{pmatrix}
1&x+1\\
0&1
\end{pmatrix}
=F^\prime_1(x,y)$ is a Jacobi matrix modulo 2 of $F$; here is how we calculate
partial derivatives modulo $2$: for instance,  
$\frac{\partial_1 g(x,y)}{\partial_1 x}=\frac{\partial_1 (y +3 x^3)}{\partial_1 x}
\cdot \frac{\partial_1 (u\XOR x)}{\partial_1 u}\big|_{u=y +3 x^3}+
\frac{\partial_1 x}{\partial_1 x}\cdot 
\frac{\partial_1 (u\XOR x)}{\partial_1 x}\big|_{u=y +3 x^3}=9x^2\cdot 1+1\cdot
1\equiv x+1\pmod 2$.
Note that a partial derivative modulo 2 of the function 
$2(x \AND y )$ is always $0$ modulo 2 because of the multiplier 2:
the function $x \AND y$ is not differentiable modulo 2 as bivariate function,
yet $2(x \AND y )$ is. So the Jacobian of the function $F$ is 
$\det F^\prime_1\equiv 1\pmod 2$.
\end{exmps*}
%
Now let  $F=(f_{1},\ldots , f_{m})\colon{\mathbb Z}^{n}_{p}\rightarrow {\mathbb Z}^{m}_{p}$  
and $f\colon{\mathbb Z}^{n}_{p}\rightarrow {\mathbb Z}_{p}$ be compatible
functions that are uniformly differentiable on $\mathbb Z_p^{n}$  modulo $p$. This is a
relatively
weak restriction since all uniformly differentiable on $\mathbb Z_p^{n}$ functions,
as well as functions that are uniformly differentiable on $\mathbb Z_p^{n}$
modulo $p^k$ for some $k\ge
1$, are uniformly differentiable on $\mathbb Z_p^{n}$ modulo $p$;
note that 
$\frac{\partial F}{\partial x_i}\equiv \frac{\partial_k F}{\partial_k x_i}\equiv
\frac{\partial_{k-1} F}{\partial_{k-1} x_i}\pmod{p^{k-1}}$.  
Moreover,
all values of all partial derivatives modulo $p^k$ (and thus, modulo $p$)
of $F$ and $f$ are $p$-adic integers everywhere on 
$\mathbb Z_p^{n}$ (see proposition \ref{intDer}),
so to calculate these values one can use the techniques considered above.

\begin{theorem}
\label{equi}
\textup{\cite{anashin2,anashin1,anashin3}}
A function  $F\colon{\mathbb Z}^{n}_{p}\rightarrow {\mathbb Z}^{m}_{p}$ is 
measure preserving whenever it is balanced modulo $p^{k}$ for some
$k\ge N_{1}(F)$ and the rank of its Jacobi matrix $F_1^\prime (\mathbf
u)$ modulo
$p$ is exactly $m$ at all points  
$\mathbf u=(u_{1},\ldots  ,u_{n})\in (\mathbb Z/p^{k}\Z)^{n}$. In case
$m=n$ these conditions are also necessary, i.e., the function $F$ preserves
measure if and only if it is bijective modulo $p^{k}$ for some
$k\ge N_{1}(F)$ and $\det(F_1^\prime (\mathbf u))\not\equiv 0\pmod{p}$ for all
$\mathbf u=(u_{1},\ldots  ,u_{n})\in (\mathbb Z/p^{k}\Z)^{n}$. Moreover,
in the considered case these conditions  imply  that $F$ preserves measure
if and only if it is bijective modulo $p^{N_1(F)+1}$. 
\end{theorem}
That is, if the mapping
$\mathbf u\mapsto F(\mathbf u)\bmod p^{N_1(F)}$ is balanced, and if
the rank of the Jacobi matrix $F_1^\prime (u)$ modulo
$p$ is exactly $m$ at all points  $
\mathbf u\in (\mathbb Z/p^{N_1(F)}\Z)^{n}$
then 
{\it each} mapping $\mathbf u\mapsto F(\mathbf u)\bmod p^r$ of
$(\mathbb Z/p^r\Z)^{n}$ onto $(\mathbb Z/p^r\Z)^{m}$
$(r=1,2,3,\ldots)$ is balanced (i.e., each point 
$\mathbf u\in (\mathbb Z/p^{r}\Z)^{m}$ has the same number of preimages
in $(\mathbb Z/p^{r}\Z)^{m}$, see definition \ref{def:erg}).
\begin{example} 
\label{KlSh-ex}
We consider as examples some mappings that were studied in \cite{KlSh} to
demonstrate how the techniques presented above work. 
\begin{enumerate}
\item {\it A mapping 
\begin{multline*}
(x,y ) \mapsto F(x,y)=\\
(x \XOR 2(x \AND y ),(y +3 x^3 )\XOR x )\bmod{2^r}
\end{multline*}
of $\mathbb (Z/2^r\Z)^{2}$ onto $\mathbb (Z/2^r\Z)^{2}$ 
is bijective for all $r=1,2,\ldots$}

Indeed, the function $F$ is bijective modulo $2^{N_1(F)}=2$ (direct verification)
and  $\det(F_1^\prime (\mathbf u))\equiv 1\pmod 2$ for all $\mathbf u\in(\mathbb
Z/2\Z)^{2}$ (see the table of derivatives in example \ref{DerLog} and examples thereafter).
\item {\it The following mappings of $\mathbb Z/2^r\Z$ onto $\mathbb Z/2^r\Z$ 
are bijective for all $r=1,2,\ldots$}: 
\begin{align*} 
x&\mapsto (x +2x^2)\bmod{2^r},\\ 
x&\mapsto (x +(x^2\OR 1))\bmod{2^r},\\ 
x&\mapsto (x \XOR (x^2\OR 1))\bmod{2^r}
\end{align*}

Indeed, all three mappings are uniformly differentiable
modulo 2, and $N_1=1$ for all of them. So it suffices to prove that
all three mappings are bijective modulo 2, i.e., as mappings of the residue
ring $\mathbb Z/2\Z$ modulo 2 onto itself (this could be checked by direct calculations), 
and that
their derivatives modulo 2 vanish at no point of $\mathbb Z/2$. The latter
also holds, since  the derivatives are, respectively,
\begin{align*}
1+4x&\equiv 1\pmod 2,\\
1+2x\cdot 1&\equiv 1\pmod 2,\\
1+2x\cdot 1&\equiv 1\pmod 2,
\end{align*}
since $(x^2\OR 1)^\prime=2x\cdot 1\equiv 1\pmod 2$, and $(x\XOR C)^\prime_1\equiv
1\pmod 2$,
see example \ref{DerLog}.
\item {\it The following closely related variants of the previous mappings
of 
$\mathbb Z/2^r$ onto $\mathbb Z/2^r$ 
are not bijective for all $r=1,2,\ldots$}:
\begin{align*}
x&\mapsto (x +x^2)\bmod{2^r},\\  
x&\mapsto (x +(x^2\AND 1))\bmod{2^r},\\ 
x&\mapsto (x +(x^3\OR 1))\bmod{2^r}, 
\end{align*}
since they are compatible but
not bijective modulo 2.
\item (see \cite{Riv}, also \cite[Theorem 1]{KlSh}) {\it Let $P (x )=a_0 +a_1 x + \cdots+a_d x^d$ be a polynomial with integral
coefficients. Then $P (x )$ is a permutation polynomial } (i.e., is bijective)
{\it modulo $2^ n$,
$n>1$ if and
only if $a_1$ is odd, $(a_2 +a_4 + \cdots)$ is even, and $(a_3 +a_5 +\cdots)$
is even.}

In view of theorem \ref{equi} we need to verify whether the two conditions
hold: first, whether $P$ is bijective modulo 2, and second,
whether
$P^\prime(z)\equiv 1\pmod 2$ for $z\in\{0,1\}$.
The first condition gives that $P(0)=a_0$ and $P(1)=a_0+a_1+a_2+\cdots a_d$
must be distinct modulo 2; hence $a_1+a_2+\cdots a_d\equiv 1\pmod 2$. 
The second condition implies that
$P^\prime(0)=a_1\equiv 1\pmod2,\ P^\prime(1)\equiv a_1+a_3+a_5+\cdots\equiv 1\pmod 2$.
Now combining all this together we get $a_2+a_3+\cdots a_d\equiv 0\pmod 2$ and 
$a_3+a_5+\cdots\equiv 0\pmod 2$, hence $a_2 +a_4 + \cdots\equiv 0\pmod 2$.
\item As a bonus, we can use exactly the same proof to
get exactly the same characterization of bijective modulo $2^r$ $(r=1,2,\ldots)$
mappings of the form $x\mapsto P (x )=
a_0\XOR  a_1x\XOR \cdots\XOR  a_dx^d\bmod 2^r$ since $u\XOR v$ is uniformly
differentiable modulo 2 as a bivariate function, and its derivative modulo
2 is exactly the same as the derivative of $u+v$, and besides, $u\XOR v\equiv
u+v\pmod 2$. 
\end{enumerate}
\end{example}
Note that in general theorem \ref{equi} could be applied to a class of
functions that is narrower than the class of all compatible functions.
However, it turns out that for $p=2$ this is not the case. Namely, the
following proposition holds, which in fact is just a restatement of a 
corresponding assertion of theorem \ref{ergBool}.
\begin{proposition}
\label{mpDer}
\textup{\cite{anashin2,anashin1}}
If a compatible function $g\colon\mathbb Z_2\rightarrow\mathbb Z_2$ preserves
measure then it is  uniformly differentiable modulo $2$ and has integer derivative
modulo $2$, which is always $1$ modulo $2$.
\end{proposition}

The techniques introduced above could also be applied to characterize ergodic
functions.
\begin{theorem}
\label{ergDer}
\textup{\cite{anashin2,anashin1,anashin3}}
Let a compatible function $f\colon{\mathbb Z}_{p}\rightarrow {\mathbb Z}_{p}$ 
be uniformly differentiable modulo $p^{2}$. Then $f$ is
ergodic if and only if it is transitive modulo $p^{N_{2}(f)+1}$ when
$p$ is an odd prime,  or modulo $2^{N_{2}(f)+2}$ when $p=2$. 
\end{theorem}
\begin{example}
\label{ergKlSh}
In \cite{KlSh} there is stated that ``...neither the invertibility nor 
the cycle structure of
$x +(x^2\OR 5)$ could be determined by his ({\slshape i.e., mine --- V.A.}) techniques.''
See however how it could be immediately done with the use of Theorem
\ref{ergDer}:
The function $f(x)=x+(x^2\OR 5)$ is uniformly differentiable
on $\mathbb Z_2$, thus, it is uniformly differentiable modulo 4 
(see example \ref{DerLog} and an example thereafter), and $N_2(f)=3$. Now to
prove that $f$ is ergodic, in view of theorem \ref{ergDer} it suffices 
to demonstrate that $f$ induces a permutation
with a single cycle on $\mathbb Z/32\Z$. Direct calculations show that a
string
$0,f(0)\bmod 32, f^2(0)\bmod 32=f(f(0))\bmod 32, \ldots, f^{31}(0)\bmod
32$ is a permutation of a string $0,1,2,\ldots,31$, thus ending the proof.
\end{example}


\section{Two fast generators}
\label{Sec:Gen}
In subsection \ref{ssec:Interp} we described how to use  interpolation series
to verify whether a given transformation $f$ of the space $\Z_2$ is ergodic (or
preserves measure): one must represent $f$ as interpolation series and apply
theorem \ref{ergBin}. Generally speaking, it is not an easy task to represent an arbitrary
continuous transformation $f$ as interpolation series (although such representation
always exists). Nevertheless, the technique works. 
Here we apply this technique 
to establish
ergodicity/measure preservation conditions for
two  special transformations that are used in cryptographic pseudorandom
generators. Both these generators are fast: The first of them uses only additions,
$\XOR$'s and multiplications by constants, the second uses additions of entries
of a certain look-up table in  accordance with bits of a variable.

\begin{theorem}
\label{thm:4.12}
The following is true:

{ $1^{\circ }$ The function $f\colon{\mathbb Z}_{2}\rightarrow {\mathbb Z}_{2}$ of the form 
$$
f(x)=a+\sum^{n}_{i=1}a_{i}(x\XOR b_{i}),
$$
\noindent where $a,a_{i},b_{i}\in {\Z}_{2}, \ \  i=1,2,3,\ldots  $, preserves measure (resp., is ergodic) if and
only if it is bijective (resp., transitive) modulo 2 (resp., modulo 4).
\par
 $2^{\circ } $ The function $f\colon{\mathbb Z}_{2}\rightarrow {\mathbb Z}_{2}$ of the form
$$
f(x)=a+\sum^{\infty }_{i=0}a_{i}\delta _{i}(x),
$$
\noindent where $a,a_{i}\in {\mathbb Z}_{2}, i=0,1,2,\ldots  $, is compatible and ergodic if and only if the
following conditions hold simultaneously:
\begin{align*}
a&\equiv 1\pmod{2};\\
a_{0}&\equiv 1\pmod{4};\\
\| a_{i}&\| _{2}=2^{-i},
\end{align*}
for $i=1,2,3,\ldots$.
{The function $f$ is compatible and measure preserving  if and only if }}
$$
\| a_{i}\| _{2}=2^{-i}\ \  (i=0,1,2,3,\ldots  ).
$$
\end{theorem}
\begin{proof}[Proof of theorem \ref{thm:4.12}] Consider interpolation series for $\delta _{i}(x)$,  $i=0,1,2,\ldots  $:
$$
\delta _{i}(x)=\sum^{\infty }_{i=0}\sigma _{i}(j){ {x}\choose {j}}.
$$
To apply theorem \ref{ergBin} we must estimate  norms of coefficients $\sigma _{i}(j)$ first. To do this, we need several lemmas.
\begin{lemma} 
\label{le:4.13} 
{ For all $i,j=1,2,3,\ldots  $ the following equations hold }
$$
\sigma _{i}(0)=0;
$$
$$
\sigma _{0}(j)=(-1)^{j+1}2^{j-1};
$$
$$
\sigma _{i}(j)=(-1)^{j+1}\sum^{\infty }_{k=1}(-1)^{k}{ {j-1}\choose {k2^{i}-1}}.
$$
\end{lemma}
\begin{proof}[Proof of lemma \ref{le:4.13}] As $\delta _{i}(0)=0$ for all $i=0,1,2,\ldots   $, then $\sigma _{i}(0)=0$.
For all $k=0,1,2,\ldots  $ we have:
 $$\delta _{i}(k)=\sum^{k}_{j=0}\sigma _{i}(j){ {k}\choose {j}}.$$
From here, with the use of formulae which express a coefficient of interpolation series of a $p$-adic function
via the values of this function in rational integer points (see e.g. \cite[Chapter
9, Section 2]{Mah}), we obtain that
$$
\sigma _{i}(j)=(-1)^{j}\sum^{\infty }_{k=0}(-1)^{k}\delta _{i}(k){ {j}\choose {k}}.
$$
{Hence, in view of the definition of the function $\delta _{i}(j)$, }
$$
\sigma _{i}(j)=(-1)^{j}\sum^{\infty }_{s=1}\sum^{s2^{i+1}-1}_{k=(2s-1)2^{i}}(-1)^{k}{ {j}\choose {k}}.
$$
{From here, using the well-known  identity (which can be easily proved) }
\begin{equation}
\label{idno1}
\sum^{n}_{k=m}(-1)^{k}{ {a}\choose {k}}=(-1)^{m}{ {a-1}\choose {m-1}}+(-1)^{n}{ {a-1}\choose {n}},
\end{equation}
{we conclude that }
$$
\sigma _{i}(j)=(-1)^{j}\sum^{\infty }_{s=1}\biggl({ {j-1}\choose {(2s-1)2^{i}-1}}-{ {j-1}\choose {2s\cdot 2^{i}-1}}\biggr).
$$
{This proves the lemma since the latter identity implies: }
$$
\sigma _{i}(j)=
\begin{cases} (-1)^{j+1}2^{j-1},&\text{if $i=0$;}\\ 
(-1)^{j+1}\sum^{\infty }_{k=1}(-1)^{k}{ {j-1}\choose {k2^{i}-1}}&\text{otherwise.} \end{cases}
$$
\end{proof}
\begin{lemma}
\label{le:4.14}
{ For all $m,t,r=0,1,2,\ldots  $ that satisfy simultaneously two
conditions $0\le t\le 2^{m}-1$ and $m\ge r$ the following congruence holds:
$$
{ {2^{m}-1}\choose {t}}\equiv (-1)^{t-\left\lfloor{t2^{-r}}\right\rfloor}{ {2^{m-r}-1}\choose {\left\lfloor{t2^{-r}}\right\rfloor}}\pmod{2^{m-r+1}}.
$$
\par In particular, for all $m,s,j\in {\mathbb N}$ that satisfy simultaneously two conditions
$m>s\ge 1$ and $j\le 2^{m-s}-1$ the following congruence holds:}
$$
{ {2^{m}-2}\choose {2^{s}j-1}}\equiv (-1)^{j}2^{s}j{ {2^{m-s}-1}\choose {j-1}}\pmod{2^{m-s+1}}.
$$
\end{lemma}
\begin{proof}[Proof of lemma \ref{le:4.14}] 
 Firstly, we recall that every $s\in {\mathbb Z}_{2}$ has a unique
representation of the form $s=2^{\ord_{2}\hskip 1pt s}\hat{s}$, where $\hat{s}$ is the unit of ${\mathbb Z}_{2}$ (i.e., $\hat{s}$
is odd, meaning $\delta _{0}(\hat{s})=1$) and henceforth has a multiplicative inverse $\hat{s}^{-1}$ in ${\mathbb Z}_{2}$. In
these denotations, assuming $M=\{i:i=1,2,\ldots  ,t;\  \ord_{2}\hskip 2pt i\ge r\}$ and $M^\prime$  a
complement of $M$ to $\{1,2,\ldots  ,t\}$, we obtain that 
\begin{multline*}
\quad{ 
{2^{m}-1}\choose {t}}=\prod^{t}_{i=1}{\frac {2^{m}-i} {i}}=\prod^{t}_{i=1}\biggl({\frac {2^{m-\ord_{2}\hskip 1pt i}} {\hat\imath}}-1\biggr)\hfill\hfill{}\equiv \\
(-1)^{\mid M^\prime \mid }\prod_{i\in M}\biggl({\hat{s}}^{-1}2^{m-\ord_{2}\hskip 1pt i}-1\biggr)\pmod {2^{m-r+1}}.
\end{multline*}
The condition $\ord_{2}i\ge r$ for $i=1,2,\ldots  ,t$ holds if and only if $i=j2^{r}$ for
$j=1,2,\ldots  ,\left\lfloor{2^{-r}t}\right\rfloor$. This means that $\mid M^\prime \mid =t-\left\lfloor{2^{-r}t}\right\rfloor$. So, the product in the
right hand part of the congruence mentioned above is equal to
\begin{multline*}
(-1)^{\mid M^\prime \mid }\prod^ {\left\lfloor{2^{-r}t}\right\rfloor}_{j=1}\bigl({\hat\jmath}^{-1}2^{m-r-\ord_{2}\hskip 1pt j}-1\bigr)=\\
(-1)^{t-\left\lfloor{t2^{-r}}\right\rfloor}{ {2^{m-r}-1}\choose {\left\lfloor{t2^{-r}}\right\rfloor}}.
\end{multline*}
This proves the first part of the statement. The second part now becomes
obvious, since
\begin{multline*}
{ {2^{m}-2}\choose {2^{s}j-1}}={\frac {2^{m}-2^{s}j} {2^{m}-1}}{ {2^{m}-1}\choose {2^{s}j-1}}\equiv \\
2^{s}j{ {2^{m}-1}\choose {2^{s}j-1}}\pmod{2^{m-s+1}}.
\end{multline*}
\end{proof}
\begin{lemma} 
\label{le:4.15} 
For $s,k=1,2,3,\ldots  $, the following holds:
\begin{enumerate}
\item $\| \sigma _{s}(k)\| _{2}\le 2^{-\left\lfloor{\log_{2}k}\right\rfloor+s-1}$,
if $k\neq 2^{s},2^{s+1}$;
\item $\|\sigma _{s}(2^{s})\| _{2}=1,\ \| \sigma _{s}(2^{s+1})\| _{2}={\dfrac {1} {2}}$;
\item $\| \sigma _{s}(2^{m}-1)\| _{2}\le 2^{-m+s-1}$, 
if $m>s\ge 1$.
\end{enumerate}
\end{lemma}
\begin{proof}[Proof of lemma \ref{le:4.15}] Represent $k$ as $k=2^{m}+t$, where $m=\left\lfloor{\log_{2}k}\right\rfloor,\  0\le t<2^{m}$.
We may assume that $m\ge s$ since otherwise $\sigma _{s}(k)=0$ in view of
lemma \ref{le:4.13}. Further, lemma \ref{le:4.13}
implies that
\begin{equation}
\label{eqno1}
\sigma _{s}(2^{m}+t)=(-1)^{t+1}\sum^{\infty }_{j=1}(-1)^{j}{ {2^{m}+t-1}\choose {2^{s}j-1}}.
\end{equation}
{With the use of the well-known identity (which can be easily  proved) }
$$
\sum^{n}_{k=0}{ {a}\choose {k}}{ {b}\choose {n-k}}={ {a+b}\choose {n}},
$$
{we obtain that}
\begin{multline}
\label{eqno2}
{ 
{2^{m}-1+t}\choose {2^{s}j-1}}=\sum^{\infty }_{k=0}{ {t}\choose {k}}{ {2^{m}-1}\choose {2^{s}j-k-1}}=\\
\sum^{\infty }_{n=0}\sum^{2^{s}-1}_{r=0}{ {t}\choose {2^{s}n+r}}{ {2^{m}-1}\choose {2^{s}(j-n-1)+(2^{s}-r-1)}}.
\end{multline}
\noindent Here, as usual, we assume that ${a \choose b}=0$ for $b<0$. In view of lemma \ref{le:4.14}, equation \eqref{eqno2} 
implies that
\begin{multline}
\label{eqno3}
\sum^\infty_{n=0}\sum^{2^{s}-1}_{r=0}(-1)^{n+r+j}{t\choose {2^sn+r}} {{2^{m-s}-1}\choose {j-n-1}}\equiv\\
{{2^{m}-1+t}\choose {2^{s}j-1}}\pmod {2^{m-s+1}}
\end{multline}
{Now \eqref{eqno1} in view of \eqref{eqno3} implies that }
\begin{multline}
\label{eqno4}
\sigma _{s}(2^{m}+t)\equiv\\ 
(-1)^{t+1}\sum^{\infty }_{n=0}\sum^{2^{s}-1}_{r=0}(-1)^{n+r}{ {t}\choose {2^{s}n+r}}\sum^{\infty }_{j=1}{ {2^{m-s}-1}\choose {j-n-1}}\hfill
\hfill{}\equiv\\ 
2^{2^{m-s}-1}(-1)^{t+1}\times\\
\sum^{\infty }_{n=0}\sum^{2^{s}-1}_{r=0}(-1)^{n+r}{ {t}\choose {2^{s}n+r}}\pmod{2^{m-s+1}}.
\end{multline}
Now applying identity \eqref{idno1} 
and assuming that $t\neq 0$, in view of lemma \ref{le:4.13} we conclude that
\begin{multline*}
(-1)^{t+1}\sum^{\infty }_{n=0}\sum^{2^{s}-1}_{r=0}(-1)^{n+r}{ {t}\choose {2^{s}n+r}}=
(-1)^{t+1}\times\\
\sum^{\infty }_{n=0}(-1)^{n}{ {t}\choose {2^{s}n+r}}
\biggl({ {t-1}\choose {2^{s}n-1}}-{ {^{t-1}}\choose {2^{s}(n+1)-1}}\biggr)=\\
2(-1)^{t+1}\sum^{\infty }_{n=1}(-1)^{n}{ {t-1}\choose {2^{s}n-1}}=2\sigma_s (t).
\end{multline*}
\noindent The left hand part of this equation is equal to -1 when $t=0$. So, taking
all these arguments into account, from \eqref{eqno4} we conclude that
$$
\sigma _{s}(2^{m}+t)\equiv 
\begin{cases} 
2^{2^{m-s}}\sigma _{s}(t)\pmod{2^{m-s+1}},&\text{if $t\neq 0$;}\\ 
-2^{2^{m-s_{-1}}}\pmod{2^{m-s+1}},&\text{if $t=0$.}
\end{cases}
$$
The latter proves statements (i) and (ii) since it easily implies
that
$$
\sigma _{s}(2^{m}+t)\equiv 
\begin{cases} 
1\pmod{2},&\text{if $m=s$, $t=0$;}\\ 
2\pmod {4},&\text{if $m=s+1$, $t=0$;}\\ 
0\pmod{2^{m-s+1}},&\text{in all other cases.}
\end{cases}
$$
Finally, if $m>s\ge 1$, then combining together lemmas \ref{le:4.13} and \ref{le:4.14}, we obtain
that
$$
\sigma _{s}(2^{m}-1)\equiv 2^{s}\sum^{\infty }_{j=1}{ {2^{m-s}-1}\choose {j-1}}\pmod{2^{m-s+1}}.
$$
{Now, applying  a well-known identity $\sum^{n}_{k=1}k{ {n}\choose {k}} =2^{n-1}n$, we conclude that }
$$
\sigma _{s}(2^{m}-1)\equiv 2^{2^{m-s}-1+s}(2^{m-s}-1)\pmod{2^{m-s+1}}.
$$
This proves (iii) and the lemma.
\end{proof}
Now everything is ready to prove  theorem \ref{thm:4.12}. We start with
the statement $1^{\circ }$. 

The operation $\XOR$ and, consequently, the function $f$ are
compatible. Now, acting as in 
we conclude that
$$
f(x)=a+\sum^{n}_{i=1}a_{i}b_{i}+\sum^{n}_{i=1}a_{i}x-2\sum^{n}_{i=1}\sum^{\infty }_{k=0}2^{k}\delta _{k}(x)\delta _{k}(b_{i}).
$$
Now, considering interpolation series for $\delta _{k}(x)$ and taking into the
account that (in view of lemma \ref{le:4.13}) $\sigma _{0}(1)=1$ and $\sigma _{i}(1)=0$ for $i=1,2,3,\ldots  $ , we obtain:
\begin{multline*}
f(x)=\\
a+\sum^{n}_{i=1}a_{i}b_{i}+ x\biggl(\sum^{n}_{i=1}a_{i}-2\sum^{n}_{i=1}\delta _{0}(b_{i})\biggr)-\sum^{\infty }_{j=2}{ {x}\choose {j}}S_{j},
\end{multline*}
where 
$
S_{j}=\sum^{n}_{i=1}\sum^{\infty }_{k=0}2^{k+1}\sigma _{k}(j)\delta _{k}(b_{i}).
$
Lemma \ref{le:4.15} immediately implies that for $k\ge 2$ 
$$
2^{k+1}\sigma _{k}(j)\equiv 
\begin{cases} 
0\pmod {2^{\left\lfloor{\log_{2}j}\right\rfloor+1}},&\text{if $j=2^{k},2^{k+1}$;}\\ 0\pmod{2^{\left\lfloor{\log_{2}j}\right\rfloor+2}},&\text{otherwise.}
\end{cases}
$$
Now theorem \ref{ergBin} implies that $f$ preserves measure (resp., is
ergodic) if and only if $\sum^{n}_{i=1}a_{i}\equiv 1\pmod {2}$ (resp., if and only if $a+\sum^{n}_{i=1}a_{i}b_{i}\equiv 1\pmod{2}$ and
$\sum^{n}_{i=1}a_{i}+2\sum^{n}_{i=1}b_{i}\equiv 1\pmod{4}$). This is obviously equivalent to the statement $1^{\circ }$
of theorem \ref{thm:4.12}.
\par
To prove statement $2^{\circ }$ of the theorem we first note that the functions $\delta _{i}$ for $i>0$ are not
compatible. As $\sigma _{i}(0)=0$ for $i>0$ (see lemma \ref{le:4.13}), we have
$$
f(x)=a+ \sum^{\infty }_{j=1}{ {x}\choose {j}}\sum^{\infty }_{i=0}a_{i}\sigma _{i}(j).
$$
Theorem \ref{ergBin} implies now that the function $f$ preserves measure if and only
if the following congruences hold simultaneously:
\begin{equation}
\label{eqnum{1}}
\begin{cases} 
\sum^{\infty }_{i=0}a_{i}\sigma _{i}(1)\equiv 1\pmod {2};&\\ 
\sum^{\infty }_{i=0}a_{i}\sigma _{i}(j)\equiv 0\pmod {2^{\left\lfloor{\log_{2}j}\right\rfloor+1}},&\text{$j=2,3,\ldots$}
\end{cases}
\end{equation}
In view of lemma \ref{le:4.13}, the first of the conditions of \eqref{eqnum{1}} is equivalent to the
congruence
\begin{equation}
\label{eqnum{2}}
a_{0}\equiv 1\pmod{2}.
\end{equation}
Moreover, lemma \ref{le:4.13} implies that $\sigma _{i}(j)=0$ for $i\ge \left\lfloor{\log_{2}j}\right\rfloor$. Hence, the
second of the conditions \eqref{eqnum{1}} is equivalent to the following system of
congruences:
\begin{equation}
\label{eqnum{3}}
\sum^{\left\lfloor{\log_{2}j}\right\rfloor}_{i=0}a_{i}\sigma _{i}(j)\equiv 0\pmod{2^{\left\lfloor{\log_{2}j}\right\rfloor+1}},\ \  j=2,3,\ldots   .
\end{equation}
Consider the following subsystem of system \eqref{eqnum{3}} for $j=2^{k}$,
$k=1,2,3,\ldots  $:
\begin{equation}
\label{eqnum{4}}
\sum^{k}_{i=0}a_{i}\sigma _{i}(2^{k})\equiv 0\pmod{2^{k+1}}, \ \ k=1,2,3,\ldots   
\end{equation}
We assert that 2-adic integers $a_{i}$ 
satisfy system of
congruences
\eqref{eqnum{4}} if and only if
$a_{i}\equiv 2^{i}\pmod{2^{i+1}}$,   $i=0,1,2,\ldots  $.
We proceed with induction on $i$. If $i=1$, then applying lemma \ref{le:4.13} for $k=1$ we conclude that
\begin{equation}
\label{eqnum{5}}
2a_{0}+a_{1}\sigma _{1}(2)\equiv 0\pmod{4}.
\end{equation}
In view of (ii) of lemma \ref{le:4.15}, the 2-adic integer $\sigma _{1}(2)$ has a multiplicative
inverse in ${\mathbb Z}_{2}$, so in view of \eqref{eqnum{2}} congruence \eqref{eqnum{5}} is equivalent to the
congruence
$$
a_{1}\equiv 2\pmod{4}.
$$
Now let the statement under the proof be true for $k<n$; consider
the congruence
\begin{equation}
\label{eqnum{6}}
\sum^{n}_{i=0}a_{i}\sigma _{i}(2^{n})\equiv 0\pmod{2^{n+1}}.
\end{equation}
By induction hypothesis, $a_{i}=2^{i}+s_{i}2^{i+1}$ $(i=0,1,\ldots  ,n-1)$ for suitable
$s_{i}\in {\mathbb Z}_{2}$. Then, taking into the account statement (ii)
of lemma \ref{le:4.15}, we conclude that
$a_{i}\sigma _{i}(2^{n})\equiv 2^{n+1}\pmod{2^{n+2}}$ for $i=0,1,\ldots  ,n-2$ and $a_{n-1}\sigma _{n-1}(2^{n})\equiv 2^{n}\pmod {2^{n+1}}$.
Hence, congruence \eqref{eqnum{6}} is equivalent to the congruence $2^{n}+a_{n}\sigma _{n}(2^{n})\equiv 0\pmod{2^{n+1}}$. As
$\sigma _{n}(2^{n})$ is a unit of ${\mathbb Z}_{2}$ (by virtue of  (ii) of
lemma \ref{le:4.15}), the latter congruence implies
that $a_{n}\equiv 2^{n}\pmod{2^{n+1}}$.
\par
From  (i) of lemma \ref{le:4.15} we easily conclude that if $a_{i}\equiv 2^{i}\pmod{2^{i+1}}$, then $a_{i}$
also satisfy each congruence of the system \eqref{eqnum{3}} for those $j$ which are not
powers of 2. This means that the set of conditions \eqref{eqnum{1}} is equivalent to
the following set of congruences:
$$
a_{i}\equiv 2^{i}\pmod{2^{i+1}},\ \  i=0,1,2,3,\ldots  .
$$
\par Thus we have proved the second part of the statement $2^{\circ }$. 
To prove the first part of this statement we note that since
$\left\lfloor{\log_{2}(i+1)}\right\rfloor+1=\left\lfloor{\log_{2}i}\right\rfloor +1$ for $i\neq 2^{k}-1$, the sufficient and necessary
conditions for the function $f$ to be ergodic (see theorem \ref{ergBin}) in the case
under consideration have the following form:
\begin{equation}
\label{eqnum{7}}
a\equiv 1\pmod{2};
\end{equation}
\begin{equation}
\label{eqnum{8}}
\sum^{\infty }_{i=0}a_{i}\sigma _{i}(1)\equiv 0\pmod{4};
\end{equation}
\begin{equation}
\label{eqnum{9}}
\sum^{\infty }_{i=0}a_{i}\sigma _{i}(j)\equiv 0\pmod{2^{\left\lfloor{\log_{2}j}\right\rfloor+1}}, \ \ j=2,3,4,\ldots   ;
\end{equation}
\begin{equation}
\label{eqnum{10}}
\sum^{\infty }_{i=0}a_{i}\sigma _{i}(2^{k}-1)\equiv 0\pmod{2^{k+1}}, \ \ k=2,3,4,\ldots .
\end{equation}
As $\sigma _{i}(1)=0$ for $i\neq 0$ (see lemma \ref{le:4.13}), then \eqref{eqnum{8}} is equivalent to the following
condition:
\begin{equation}
\label{eqnum{11}}
a_{0}\equiv 1\pmod{2}.
\end{equation}
During the proof of the second part of the statement $2^{\circ }$ we have
established that if $a_{0}\equiv 1\pmod{2}$ (and, in particular, if \eqref{eqnum{11}} is
satisfied) then the conditions \eqref{eqnum{9}} are equivalent 
to conditions
\begin{equation}
\label{eqnum{12}}
a_{i}\equiv 2^{i}\pmod{2^{i+1}},\ \  (i=1,2,3,\ldots  ).
\end{equation}
Finally, combining together statements  (i) of lemma \ref{le:4.15} and of
lemma \ref{le:4.13} we conclude that that if 2-adic integers $a_{i}$
$(i=0,1,2,\ldots  )$ satisfy conditions \eqref{eqnum{12}} and \eqref{eqnum{11}} simultaneously, then
$a_{i}$ also satisfy conditions \eqref{eqnum{10}}. Thus, the  union  of  conditions  \eqref{eqnum{7}}---\eqref{eqnum{10}}  is equivalent to the union of conditions \eqref{eqnum{7}}, \eqref{eqnum{11}}, and of \eqref{eqnum{12}}. This
proves the first part of the statement $2^{\circ }$ and the whole theorem
\ref{thm:4.12}.
\end{proof}

\section{Estimates of randomness}
\label{sec:Rand}

Loosely speaking, within a context of this paper a PRNG is an algorithm that
takes a short binary word (an initial state, a seed) and stretches it to
a much longer word, which for any seed must look like random, that is, like
a sequence of fair coin tosses. Given a seed,
the whole period of
the produced sequence (which is necessarily periodic) is never used in practice.
However, the period must be very long and as `random-looking' as possible.
In most applications (e.g., in cryptography), a period of the output sequence much be exponentially longer than the seed, and the algorithm must be fast; whence, the corresponding program
cannot be complicated. Thus, designing a PRNG is a kind of paradox: On the one hand, the outputted string must `look like random' (say, must have
high Kolmogorov complexity); on the other hand, the generating program
must be short, whence, the Kolmogorov complexity of the produced sequence
will be necessarily low.

In real life settings they often agree that the  output sequence `looks sufficiently
random' whenever it
passes certain (in some cases, rather limited) number of statistical tests. In particular,
the outputted  string must have no obvious structure using which one can,  given a segment of the output sequence, 
 predict
with high probability the next bit.
Of course, at least
{\it some} sequences generated by compatible ergodic transformations of the
space $\mathbb Z_2$ are highly predictable, e.g., sequences (even truncated
ones) produced by linear congruential generators,  see \cite{Shn} and references
therein. Note that recently there were developed a number of effective prediction methods
for machine learning,  e.g. transduction
\cite{Vap}, conformal prediction and some others, see \cite{Vovk}. It would
be very interesting to understand what sequences generated by compatible ergodic
transformations of the space $\mathbb Z_2$ can be predicted by these methods. However, this question is outside
the scope of the given paper and can be a theme of a future work.   
      
In this section we pursue a much
less ambitious goal: We study distributions and structural properties  of sequences
produced by compatible ergodic transformations of the space $\mathbb Z_2$ in order to
demonstrate that at least with respect to {\it some} tests based on distribution
of patterns these sequences are good.

A word of caution:  For some convenience during proofs, throughout this section speaking
of base-2 expansions, as well as of 2-adic representations, we read them {\it from left to right}, so $1101$ means $1101000\ldots$; and 1101 is a
base-2 expansion of
11, and {\it not} of 13!

\subsection{Distribution of $k$-tuples} 

Whenever $f$ is a compatible ergodic transformation 
of the space $\Z_2$, the sequence $\mathcal T_n=\{z_i=f^i(z)\mod 2^n\}_{i=0}^\infty$ is 
strictly uniformly distributed  as a sequence of binary words of
length $n$ (see section \ref{sec:App}). 
However, for applications it is important to study distributions of
a binary sequence $\mathcal T_n^\prime$ obtained from $\mathcal T$ by concatenation
of these $n$-bit words: 
%
However, one could consider
the same sequence as a binary sequence and ask what is a distribution
of $n$-tuples in this binary sequence. {\it Strict uniform distribution of an
arbitrary sequence $\mathcal T$
as a sequence over $\Z/2^n\Z$ does not necessarily imply uniform distribution
of overlapping $n$-tuples, if this sequence is considered as a binary sequence!} 

For instance, let $\mathcal T$ be the following strictly uniformly
distributed sequence over $\Z/4\Z$ with period length exactly $4$:
$\mathcal T=023102310231\ldots$. Then its representation as a binary sequence
is $\mathcal T^\prime=000111100001111000011110\ldots$ 
Obviously, when we consider $\mathcal T$ as
a sequence over the residue ring $\Z/4\Z$, then each number of $\{0,1,2,3\}$ occurs in the
sequence with the same frequency $\frac{1}{4}$. Yet if we consider $\mathcal
T$ as a binary sequence, then $00$ (as well as $11$) occurs in this sequence with
frequency $\frac{3}{8}$, whereas $01$ (and $10$) occurs with frequency $\frac{1}{8}$.
Thus, the sequence $\mathcal T$ is uniformly distributed over $\Z/4\Z$, and
it is not uniformly distributed over $\Z/2\Z$.

In this subsection we show that this effect does not take
place for the sequences $\mathcal T_n$:
{\it Considering this sequence as
a binary sequence, a distribution
of $k$-tuples is uniform, for all $k\le n$}. Now we state this property more formally.

Consider a (binary) {\it $n$-cycle} 
$C=(\varepsilon_0\varepsilon_1\dots \varepsilon_{n-1})$; 
that is, an oriented
graph with vertices $\{a_0,a_1,\ldots, a_{n-1}\}$ and edges 
$$\{(a_0,a_1),(a_1,a_2),\ldots, (a_{n-2},a_{n-1}),(a_{n-1},a_0)\},$$ 
where
each vertex $a_j$ is labelled with $\varepsilon_j\in\{0,1\}$, $j=0,1,\dots,n-1$.
(Note that then $(\varepsilon_0\varepsilon_1\dots \varepsilon_{n-1})=
(\varepsilon_{n-1}\varepsilon_0\dots \varepsilon_{n-2})=\ldots$, etc.).

Clear, each purely periodic sequence $\mathcal S$ over $\Z/2\Z$ with period 
$\alpha_0\ldots\alpha_{n-1}$
of length $n$
could be related to a binary $n$-cycle $C(\mathcal S)=(\alpha_0\ldots\alpha_{n-1})$.
Conversely, to each binary $n$-cycle $(\alpha_0\ldots\alpha_{n-1})$ we could
relate $n$ purely periodic binary sequences of period length $n$: They
are $n$ shifted versions of the sequence
$$\alpha_0\ldots\alpha_{n-1}\alpha_0\ldots\alpha_{n-1}\ldots,$$
that is
\begin{align*}
&\alpha_1\ldots\alpha_{n-1}\alpha_0\alpha_1\ldots\alpha_{n-1}\alpha_0\ldots,\\
&\alpha_2\ldots\alpha_{n-1}\alpha_0\alpha_1\alpha_2\ldots\alpha_{n-1}\alpha_0\alpha_1\ldots,\\
&\ldots\qquad\ldots\qquad\ldots\\
&\alpha_{n-1}\alpha_0\alpha_1\alpha_2\ldots\alpha_{n-2}\alpha_{n-1}\alpha_0\alpha_1\alpha_2\ldots\alpha_{n-2}\ldots
\end{align*}

Further, {\it a $k$-chain in a binary $n$-cycle}  
$C$ is a
binary string $\beta_0\dots\beta_{k-1}$, $k<n$, that satisfies the following
condition: There exists $j\in\{0,1,\ldots,n-1\}$ such that $\beta_i=\varepsilon_{(i+j)\bmod
n}$ for $i=0,1,\ldots, k-1$. Thus, a $k$-chain
is just a string of length
$k$ of labels that corresponds to a chain of length $k$ in a graph $C$.

We call a binary $n$-cycle $C$ {\it $k$-full}, if each $k$-chain
occurs in the graph $C$ the same number $r>0$ of times.

Clearly, if $C$ is $k$-full, then $n=2^kr$. For instance, a well-known
De Bruijn sequence is an $n$-full $2^n$-cycle. 
It is clearly that a $k$-full $n$-cycle is $(k-1)$-full:
Each $(k-1)$-chain occurs in $C$ exactly $2r$ times, etc. Thus, if an $n$-cycle
$C(\mathcal S)$ is $k$-full, then each $m$-tuple (where $1\le m\le k$) occurs in
the sequence $\mathcal S$ with the same probability (limit frequency) $\frac{1}{2^m}$.
That is, the sequence $\mathcal S$ is {\it $k$-distributed}, see
\cite[Section 3.5, Definition D]{knuth}.
\begin{definition} A purely periodic binary sequence $\mathcal S$ with period length
exactly $N$ is said to be {\it
strictly $k$-distributed} if and only if a corresponding $N$-cycle $C(\mathcal S)$
is $k$-full.
\end{definition}

Thus, if a sequence $\mathcal S$ is strictly $k$-distributed, then it is
strictly $s$-distributed, for all positive $s\le k$.

 A $k$-distribution is a good `indicator
of randomness' of an infinite sequence: The larger $k$, the better the
sequence, i.e., `more random'. The best case is when a sequence is $k$-distributed
for all $k=1,2,\ldots$. Such sequences are called $\infty$-distributed.
Obviously, a periodic sequence can not be $\infty$-distributed.

On the other hand, a periodic sequence is just an infinite repetition of a finite
sequence, the period. 
So we
are interested in `how random' this finite sequence (the period) is.
Of course, it seems very reasonable to consider a period of length $n$ as an $n$-cycle 
and to study a distribution
of $k$-tuples in $n$-cycle; for instance,
if this $n$-cycle is $k$-full, the distribution of $k$-tuples is strictly
uniform. However, other approaches also exist.

In \cite[Section 3.5, Definition Q1]{knuth} there is considered the following 
`indicator of randomness'
of a finite sequence over a finite alphabet $A$ (we formulate the corresponding
definition for $A=\{0,1\}$): 
a finite binary sequence 
$\varepsilon_0\varepsilon_1\dots \varepsilon_{N-1}$
of length $N$ is said to be random (sic!), if and only if
\begin{equation}
\label{eq:Q1}
\bigg|\frac{\nu(\beta_0\ldots\beta_{k-1})}{N}-\frac{1}{2^k}\bigg|\le\frac{1}{\sqrt
N}
\end{equation}
for all $0<k\le\log_2N$, where $\nu(\beta_0\ldots\beta_{k-1})$ is the number
of occurrences of a binary word $\beta_0\ldots\beta_{k-1}$ in a binary word
$\varepsilon_0\varepsilon_1\dots \varepsilon_{N-1}$. If a finite sequence
is random in the meaning of this Definition Q1 of \cite{knuth}, we shall say
that it has {\it a property} Q1, or {\it satisfies} Q1. We shall also
say that an {\it infinite periodic sequence satisfies} Q1 if and only if its exact
period satisfies Q1.
Note that, contrasting to the case of strict $k$-distribution, which implies
strict $(k-1)$-distribution, 
it is not enough to demonstrate only
that inequality \eqref{eq:Q1}
holds for $k=\lfloor\log_2N\rfloor$ to prove a finite sequence of length $N$ 
satisfies Q1:
For instance, a sequence $1111111100000111$ satisfies \eqref{eq:Q1} for
$k=\lfloor\log_2n\rfloor=4$, and does not satisfy \eqref{eq:Q1} for $k=3$.
Note that an analog of property Q1 for odd prime $p$ could be stated in an obvious
way. 

Now we are able to state the following theorem.
\begin{theorem}
\label{thm:distr}
Let $\mathcal T^\prime_n$ be a
binary representation of the sequence $\mathcal T_n$ {\rm (hence $\mathcal T^\prime_n$ is
a purely periodic binary sequence of period length exactly $n2^n$)}.
Then 
the sequence $\mathcal T^\prime_n$ is strictly $n$-distributed.
Moreover, 
this sequence satisfies
{\rm Q1}.
\end{theorem}
\begin{proof} 
Let
$\mathcal T^\prime_n=\zeta_0\zeta_1\ldots$ be a binary representation of
the sequence $\mathcal T_n$. Take
an arbitrary binary word $\mathbf b=\beta_0\beta_1\ldots\beta_{n-1}$, $\beta_j\in\{0,1\}$,
and for $k\in\{0,1,\ldots, n-1\}$ denote 
\begin{multline*}
\nu_k(\mathbf b)=\\
|\{r\: 0\le r<n2^n;
r\equiv k\pmod n;\\
\zeta_r\zeta_{r+1}\ldots\zeta_{r+n-1}=
\beta_0\beta_1\ldots\beta_{n-1}\}|
\end{multline*}
Obviously, $\nu_0(\mathbf b)$ is the number of occurrences of a rational
integer $z$ with base-$2$ expansion $\beta_0\beta_1\ldots\beta_{n-1}$ at
the exact period of the sequence $\mathcal Z$. Hence, $\nu_0(\mathbf b)=1$
since the sequence $\mathcal T_n$ is strictly uniformly distributed modulo
$2^n$. Now consider $\nu_k(\mathbf b)$ for $0<k<n$.

Fix $k\in\{1,2\ldots,n-1\}$ and let $r=k+tn$. Since $f$ is compatible,  then 
$\zeta_r\zeta_{r+1}\ldots\zeta_{r+n-1}=\beta_0\beta_1\ldots\beta_{n-1}$
holds if and only if the following two relations hold simultaneously:
\begin{equation}
\label{eq:distr1}
\zeta_{tn+k}\zeta_{tn+k+1}\ldots\zeta_{tn+n-1}=\beta_0\beta_1\ldots\beta_{n-k-1}
\end{equation}
\begin{multline}
\label{eq:distr2}
f_{t}(\overline{\zeta_{tn}\zeta_{tn+1}\ldots\zeta_{tn+k-1}})\equiv\\
\overline{\beta_{n-k}\beta_{n-k+1}\ldots\beta_{n-1}}\pmod{2^k}.
\end{multline}
Here $\overline{\gamma_0\gamma_1\ldots\gamma_s}=\gamma_0+\gamma_1\cdot
2+\dots+\gamma_s\cdot 2^s$ for $\gamma_0,\gamma_1,\ldots,\gamma_s\in\{0,1\}$
is a rational integer with a base-$2$ expansion $\gamma_0\gamma_1\ldots\gamma_s$.

For a given 
$\mathbf b=\beta_0\beta_1\ldots\beta_{n-1}$
congruence \eqref{eq:distr2} has exactly one solution 
$\overline{\alpha_0\alpha_1\dots\alpha_{k-1}}$ modulo $2^k$, since
$f$
is ergodic, whence, bijective modulo $2^k$.
Thus,
in view of \eqref{eq:distr1} and \eqref{eq:distr2} we conclude that 
$\zeta_r\zeta_{r+1}\ldots\zeta_{r+n-1}=\beta_0\beta_1\ldots\beta_{n-1}$
holds if and only if 
\begin{equation}
\label{eq:distr3}
\zeta_{s}\zeta_{s+1}\ldots\zeta_{s+n-1}=
\alpha_0\alpha_1\dots\alpha_{k-1}\beta_0\beta_1\ldots\beta_{n-k-1},
\end{equation}
where $s=tn$. Yet there 
exists exactly one 
$s\equiv 0\pmod n$, $0\le s< 2^nn$ such that  \eqref{eq:distr3} holds, 
since every element of
$\Z/2^n\Z$ occurs at the period of $\mathcal T_n$ exactly once. We conclude
now that
$\nu_k(\mathbf b)=1$ for all $k\in\{0,1,\ldots, n-1\}$; thus, $\nu(\mathbf b)=
\sum_{j=0}^{n-1}\nu_j(\mathbf b)=n$ for all $\mathbf b$. This means that
the $(n2^n)$-cycle 
$C(\mathcal T^{\prime}_n)$ is $n$-full, whence, the sequence $\mathcal T^{\prime}_n$
is strictly $n$-distributed.
This completes the proof of the first assertion of the theorem.

To prove the second assertion note that 
in view of the first assertion every $m$-tuple 
for $1\le m\le n$ occurs at the $n2^n$-cycle $C(\mathcal T^\prime_n)$ exactly
$2^{n-m}n$ times. Thus, every such $m$-tuple occurs $2^{n-m}n-c$ times
in the finite binary sequence 
$\hat{\mathcal T}_n=\hat z_0\hat z_1\ldots\hat z_{2^n-1}$, where
$\hat z$ for $z\in\{0,1,\ldots,2^n-1\}$ is an $n$-bit sequence that agrees
with base-$2$ expansion of $z$. Note that $c$ depends on the $m$-tuple, yet
$0\le c\le m-1$ for every $m$-tuple. Easy algebra shows that \eqref{eq:Q1}
holds for these $m$-tuples. 

Now to prove that $\mathcal T^\prime_n$ satisfies Q1 
we have only to demonstrate that \eqref{eq:Q1} holds for $m$-tuples with
$m=n+d$, where $0<d\le\log_2n$. We claim that any such $m$-tuple occurs in
the sequence $\hat{\mathcal T}_n$ not more than $n$ times.

Indeed, in this case 
$\zeta_r\zeta_{r+1}\ldots\zeta_{r+n+d-1}=\beta_0\beta_1\ldots\beta_{n+d-1}$
holds if and only if besides the two relations \eqref{eq:distr1} and \eqref{eq:distr2} the following
extra congruence holds:
\begin{multline*}
f(\overline{\zeta_{tn}\zeta_{tn+1}\ldots\zeta_{tn+k-1}\beta_0\beta_1\ldots\beta_{d-1}})
\equiv\\
\overline{\beta_{n-k}\beta_{n-k+1}\ldots\beta_{n+d-1}}\pmod{2^{k+d}},
\end{multline*}
where $k=r\bmod n$. Yet this extra congruence may or may not have a solution
in unknowns $\zeta_{tn},\zeta_{tn+1},\ldots,\zeta_{tn+k-1}$; this depends on $\beta_0\beta_1\ldots\beta_{n+d-1}$.
But if such  solution exists, it is unique for a given $k\in\{0,1,\ldots,n-1\}$, since $f$
is ergodic, whence, bijective modulo $2^s$ for all $s=1,2,\ldots$.
This proves our claim. Now  exercise in inequalities shows that \eqref{eq:Q1}
holds in this case, thus completing the proof of the theorem.
\end{proof}
\begin{note}
\label{note:distr} 

The second assertion of theorem \ref{thm:distr} holds for arbitrary prime
$p$. Namely, {\it a base-$p$ representation of an output sequence 
of a congruential generator over $\Z/p^n\Z$
of a maximum period length is strictly $n$-distributed sequence over $\Z/p\Z$
of period length exactly $p^nn$, which satisfies Q1}.

Moreover, the first assertion of theorem \ref{thm:distr} also holds 
for a \emph{truncated} congruential generator; that is, for a generator $\mathfrak
A$ of section \ref{sec:App} with output
function $F(x)=\big\lfloor\frac{x}{p^{n-k}}\big\rfloor\bmod p^k$. Namely, {\it
a base-$p$
representation of  the output sequence of a truncated congruential generator
over $\Z/p^n\Z$ of a maximum period length is a purely periodic strictly $k$-distributed
sequence over $\Z/p\Z$ of period length $p^nk$}.

The second assertion for this generator holds whenever $2+p^k>kp^{n-k}$;
thus, {\it  one could truncate $\le\big(\frac{n}{2}-\log_p\frac{n}{2}\big)$ lower order digits 
without affecting property Q1}.

All these statements could be proved by slight modifications of the
proof of theorem \ref{thm:distr}. We omit details.
\end{note}

\subsection{Coordinate sequences}

In this subsection, we study some structural properties 
of a binary sequence produced by a compatible ergodic transformation $f$
of the space
$\mathbb Z_2$. 
Clear, a binary sequence $\mathcal
S_j=\{\delta_j(f^i(z_0)\}_{i=0}^\infty$ (which is called the \emph{$j$-th
coordinate sequence}, is a purely periodic binary sequence of period length $2^{j+1}$. 

Moreover, it easy to understand that \emph{the second half of the period
of every coordinate sequence $\mathcal S_j=s_0,s_1,s_2,\dots$ is a bitwise negation of its first half}: 
\begin{equation}
\label{eq:halfper}
s_{i+2^j}\equiv s_i+1\pmod 2, \ \ i=0,1,2,\ldots
\end{equation}
This
immediately follows from theorem \ref{ergBool} and means, loosely speaking,
that the $j$-th 
coordinate sequence is as
complex as the first half of its period. So  it is important to know what
sequences of length $2^j$ could be outputted as the first half of the
period of the $j$-th 
coordinate sequence; more formally, what
values are taken  by the rational integer $\gamma=s_0+s_12+s_22^2+\dots+s_{2^{j}-1}2^{2^{j}-1}$,
for the $j$-th coordinate sequence
$\mathcal S_j=s_0,s_1,s_2,\dots$. 

In other words, let $\gamma_j(f,z)\in\mathbb N_0$ be such a number that its base-$2$ expansion 
agrees with the first half
of the period of the $j$\textsuperscript{th} coordinate sequence; 
i.e., let
\begin{multline*}
\gamma_j(f,z)=\delta_j(f^{0}(z))+2\delta_j(f^{1}(z))+
4\delta_j(f^{2}(z))+\cdots\\
+2^{2^j-1}\delta_j(f^{2^j-1}(z)).
\end{multline*}
Obviously, $0\le\gamma_j(f,z)\le 2^{2^j}-1$. The following natural question should be answered:
{\it Given a compatible and ergodic mapping
$f\colon\mathbb Z_2\rightarrow\mathbb Z_2$ and a $2$-adic integer $z\in\mathbb
Z_2$, what infinite string $\gamma_0=\gamma_0(f,z),\gamma_1=\gamma_1(f,z),
\gamma_2=\gamma_2(f,z),\dots$ (where $\gamma_j\in\{0,1,\dots,2^{2^j}-1\}$
for
$j=0,1,2,\dots$) could be obtained?} 

And the answer is:  {\it any one.} Namely, the
following theorem holds (which, interestingly, could be proved by a `purely 2-adic' argument). 

\begin{theorem}
\label{AnyHalfPer}
Let $\Gamma=\{\gamma_j\in\mathbb N_0\colon  j=0,1,2,\ldots\}$
be an arbitrary sequence of non-negative rational integers that satisfy
$0\le\gamma_j\le 2^{2^j}-1$ for $j=0,1,2,\ldots$.
There exists a compatible and ergodic mapping 
$f\colon\mathbb Z_2\rightarrow\mathbb Z_2$ and a $2$-adic integer 
$z\in\mathbb Z_2$ such that $\delta_j(z)=\delta_0(\gamma_j)$, 
$\delta_0(f^{i}(z))\equiv \gamma_0+i\pmod 2$, 
and
$$\delta_j(f^{i}(z))\equiv \delta_{i\bmod{2^j}}(\gamma_j)+
\biggl\lfloor\frac{i}{2^j}\biggr\rfloor\pmod 2$$ for all
$i,j\in\mathbb N$. 
\end{theorem}
\begin{note*} The sequence 
$\Bigl\{\Bigl\lfloor\frac{i}{2^j}\Bigr\rfloor\bmod 2\: i=1,2,
\ldots\Bigr\}$ is merely a binary sequence of alternating gaps and runs (i.e., blocks
of consecutive $0$'s or $1$'s, respectively) of length
$2^j$ each.
\end{note*}
\begin{proof}[Proof of theorem \ref{AnyHalfPer}]
Put $z=z_0=\sum_{j=0}^{\infty}\delta_0(\gamma_j)2^j$ and
\begin{multline*} 
z_i= (\gamma_0+i)\bmod 2+\\
\sum_{j=1}^{\infty}\biggl(\biggl(
\delta_{i\bmod{2^j}}(\gamma_j)+
\biggl\lfloor\frac{i}{2^j}\biggr\rfloor\biggr)\bmod 2\biggr)\cdot
2^j
\end{multline*}
for $i=1,2,3,\ldots$ . Consider a sequence $Z=\{z_i\colon i=0,1,2,\ldots\}$.
Speaking informally, we are filling a table with countable infinite number of rows
and columns in such a way that the first $2^j$ entries of the $j$-th 
column represent $\gamma_j$ in its base-2 expansion, and the other entries
of this column are obtained from these by applying recursive relation  \eqref{eq:halfper}.
Then each $i$\textsuperscript{th} row of the table is a 2-adic canonical
representation of $z_i\in Z$.

We shall prove that $Z$ is a dense subset in $\mathbb Z_2$, and then
define $f$ on $Z$ in such a way that $f$ is compatible and ergodic on $Z$.
This will imply the assertion of the theorem. 

Proceeding along this way we claim that $Z\bmod 2^k = \mathbb Z/2^k\Z$ for all $k=1,2,3,\ldots$,
i.e., a natural ring homomorphism $\bmod\, 2^k\colon z\mapsto z\bmod 2^k$ maps
$Z$ onto the residue ring $\mathbb Z/2^k\Z$. Indeed, this trivially holds
for $k=1$. Assuming our claim holds for $k< m$ we prove it for $k=m$.
Given arbitrary $t\in\{0,1,\ldots,2^{m}-1\}$ there exists $z_i\in Z$ such
that $z_i\equiv t\pmod{2^{m-1}}$. If $z_i\not\equiv t\pmod{2^{m}}$ then
$\delta_{m-1}(z_i)\equiv\delta_{m-1}(t)+1\pmod 2$ and thus
$\delta_{m-1}(z_{i+2^{m-1}})\equiv\delta_{m-1}(t)\pmod 2$. However, 
$z_{i+2^{m-1}}\equiv z_i\pmod {2^{m-1}}$. Hence
$z_{i+2^{m-1}}\equiv t\pmod {2^m}$. 

A similar argument shows that for each $k\in\mathbb N$ 
the sequence $\{z_i\bmod 2^k\}_{i=0}^\infty$
is purely periodic with period length $2^k$, and each $t\in\{0,1,\ldots,2^{k}-1\}$
occurs at the period exactly once (in particular, all members of $Z$ are
pairwise distinct 2-adic integers). Moreover, $i\equiv i^{\prime}\pmod{2^k}$
if and only if $z_{i}\equiv z_{i^{\prime}}\pmod{2^k}$. Consequently, $Z$ is dense
in $\mathbb Z_2$ since for each $t\in\mathbb Z_2$ and each $k\in\mathbb
N$ there exists $z_i\in Z$ such that $\|z_i-t\|_2\le 2^{-k}$. Moreover, if
we define $f(z_i)=z_{i+1}$ for all $i=0,1,2,\ldots$ then 
$\|f(z_i)-f(z_{i^{\prime}})\|_2=\|z_{i+1}-z_{i^{\prime}+1}\|_2=
\|(i+1)-(i^{\prime}+1)\|_2=\|i-i^{\prime}\|_2=\|z_i-z_{i^{\prime}}\|_2$.
Hence, $f$ is well defined and compatible on $Z$; it follows that the continuation
of $f$ to the whole space $\mathbb Z_2$ is compatible. Yet $f$ is transitive
modulo $2^k$ for each $k\in\mathbb N$, so its continuation is ergodic.
\end{proof}
\section{Conclusion}
\label{sec:Conc}
 
In this paper, we demonstrate that, loosely speaking,
a contemporary digital computer `thinks 2-adically': Most common processor
instructions, both numerical (i.e., arithmetic, e.g. addition, multiplication),
logical (such as bitwise $\OR$, $\AND$, $\XOR$, $\NOT$) and machine (left
and right shifts) are continuous functions with respect to 2-adic metric.
Hence, a computer program which is combined from these operators is a continuous
function defined on (and valuated in) the space of 2-adic integers. So we
believe that natural metric
for a digital computer is non-Archimedean: The sequence of states of a
program (as we have demonstrated by example of programs that generate pseudorandom numbers)
admits an adequate description as a smooth trajectory in the non-Archimedean metric space. If so, a digital computer is likely to be  perfect for simulating non-Archimedean
dynamics, and not as good for simulating Archimedean systems.

The later phenomenon was already noticed in numerical analysis:
For instance, paper
\cite{Li} reads:
\begin{quote}
Digital computers are absolutely incapable of showing true long-time
dynamics of some chaotic systems, including the tent map, the Bernoulli shift map
and their analogues, even in a high-precision floating-point arithmetic. 
\end{quote}

Note that both these dynamical systems, the tent map and the Bernoulli shift map,
are ergodic. However, theoretical analysis, as well as 1000
computer verifications in \cite{Li} demonstrate that behaviour of corresponding
computer programs is {\it not} ergodic:

\begin{quote}
It is found that all chaotic orbits will be eventually converge to zero within $N_r$ iterations, and
that the value of $N_r$ is uniquely determined by the details of digital floating-point arithmetic. 
\end{quote}

Moreover, inspired by results of \cite{Li} we undertook our own study  of  {\it discrete versions} of these two maps, supported
by
computer
experiments based on fixed-point (actually, integer) arithmetic instead of floating-point one. Namely, we considered a map
$B_n\colon x\mapsto \frac{(x \OR 1)-1}{2}\pmod {2^n}$ as a discrete analog of the Bernoulli
shift map, and a map  $T_n\colon x\mapsto \frac {x \AND (-2)}{2}-x\cdot(x\AND 1)\pmod
{2^n}$, as a discrete analog of the tent map. Both these maps are transformations
of the set $\{0,1,\ldots, 2^n-1\}=\mathbb Z/2^n \mathbb Z$, and elements of latter set can be put into a correspondence
with real numbers in $[0,1]$ via the Monna map, 
$$x=\sum_{i=0}^{n-1}\delta_i(x)2^i\longleftrightarrow
\sum_{i=0}^{n-1}\delta_i(x)2^{-i-1}
\in [0,1].$$

e.g., $2=\ldots 0010\longleftrightarrow\frac{1}{4}$, $3=\ldots0011\longleftrightarrow  \frac{1}{2}+\frac{1}{4}=\frac{3}{4}$,
etc. 
Up to this correspondence, both $B_n$ and $T_n$ give the same plots in a
unit square as, respectively, the Bernoulli shift and the tent map, being restricted
to real numbers with $n$ binary digits after the point. However, both $B_n$
and $T_n$ are {\it not} ergodic either: $B_n$ converges to 0 after at most $n$ iterations,
and $T_n$ always falls in short cycles, of length $n$ at most.

This effect cannot occur for truly ergodic maps: Loosely speaking, ergodic
transformations have no invariant subsets, except of subsets of measure 0
and of full measure. 
Thus, any ergodic transformation of a finite set (which is endowed with a natural probabilistic uniform
measure) must necessarily be transitive, i.e., must permute all elements
of the set cyclically. In other words, these considerations show that
computer simulations of Archimedean ergodic systems are indeed  inadequate, since
the corresponding programs clearly exhibit a non-ergodic behaviour.

On the contrary, results of the present paper  demonstrate that whenever one
considers ergodic transformation of the space of 2-adic integers that satisfy
Lipschitz condition with a constant 1, any restriction of this transformation
to $n$-bit precision remains ergodic: 
Thus, digital computers are
perfect for simulating behaviour of these 2-adic dynamical systems: In the
paper, the corresponding
dynamics was  used to construct effective pseudorandom generators with
prescribed characteristics. Numerous computer experiments with these programs
(e.g., the ones undertaken during the development of the ABC stream cipher
\cite{abc-v2}) are in full agreement with the theory presented above. At our view, these considerations give us another
evidence  that a non-Archimedean
(namely, 2-adic) metric is natural for digital computers, whereas the Archimedean metric
is not.

Yet another evidence is given by the following observation: Every digital computer, even the simplest one, can, by its very origin, properly operate with 2-adic
numbers.
Let's undertake the following `computer experiment'. Start MS Windows XP, run
a built-in Calculator. Switch to Scientific mode. Press {\tt Dec} (that
is, switch to decimals), press 1, then +/-.  The calculator  returns -1, as prescribed. 

Now, press {\tt
Bin}, switching the calculator to binaries. The calculator
returns ...111 (64 ones),
a 2-adic representation of -1, up to the highest precision the calculator
could achieve, 64 bits. (Here a programmer will most likely say that the calculator just uses
the {\it two's complement}).

Now press {\tt Dec} again; the calculator returns 18446744073709551615.
This number is congruent to -1 modulo $2^{64}$. Now press
successively 
{\tt / }, {\tt 3}, {\tt =}, {\tt Bin}, thus dividing the number by 3 and representing the result in a binary form.  The calculator returns ...10101010101, a 2-adic representation
of -1/3, with 2-adic precision $2^{-64}$. Indeed, switching back to {\tt
Dec} we obtain 6148914691236517205, a multiplicative inverse to -3 modulo
$2^{64}$.

This toy experiment could be performed on most calculators. However, sometimes
a calculator returns an erroneous result. This usually happens when
a corresponding program is written in a higher-order language. 
Very loosely speaking, the capability of a calculator to perform
2-adic arithmetic depends on how the corresponding program is written:
programs written in assembler usually are more capable to perform 2-adic calculations than 
the ones written in higher-level languages. Programmers use assembler
when they want to exploit CPU's resources in the most optimal way; e.g., to store
negative numbers they use
the two's complement rather than reserve special registry for a sign.  But the
usage of the
two's complement of $x$ (that is, of $\NOT x$) is just a way to represent a negative integer in a 2-adic
form, $-x=1+\NOT x$, see equations \eqref{eq:id} of Section \ref{sec:App}.  Thus, we might conclude that a
CPU is used in a  more optimal way when it actually works with binary words
as with 2-adic numbers.
Thus, a CPU looks more  `non-Archimedean-oriented'  
than `Archimedean-oriented'. 

We human beings are Archimedean creatures: We agree that the surrounding  physical world is Archimedean judging by
numerous experiments. Our experience gives us a strong evidence
that trajectories
of a physical (especially, mechanical) dynamical system admit (as a rule) adequate descriptions  by smooth curves in an Archimedean (Euclidean) metric
space.  Moreover, we can simulate behaviour of these mechanical systems
by other physical processes, e.g., by electrical ones: This way we come to
{\it analog} computers  that can simulate processes of our physical (at least,
mechanical) world with arbitrary high precision since their internal basic
operators
are continuous functions with respect to Archimedean metric.

But then,
if we see that a digital computer cannot simulate long-time dynamics even
of rather simple Archimedean dynamical systems, yet can simulate with arbitrarily
high precision non-Archimedean dynamics, we probably should agree that 
digital computers are a kind of non-Archimedean devices, something like {\it
analog computers for the non-Archimedean world}, since their internal basic
operators are continuous functions with respect to the 2-adic (i.e., non-Archimedean)
metric.


We believe that these considerations must be taken into account while simulating
dynamical systems on digital computers: Probably, the simulation will be
adequate for non-Archimedean dynamical systems, whereas for non-Archimedean
ones it will be not.

Also, the approach presented in the paper could probably be applied to other
problems of computer science, and not only to the problem of pseudorandom
generation. For instance, consider an automaton with
a binary input and binary output. This automaton actually performs
a transformation of the space $\mathbb Z_2$ of 2-adic integers: Each infinite
input string of 0s and 1s the automaton transforms into an infinite output
string of 0s and 1s (we suppose that the initial state is fixed). Note that every outputted  $i$-th bit depends only on the inputted $i$-th bit and on the current state
of the automaton. Yet the current state depends only on the previous state
and on the $(i-1)$-th input bit. Hence, for every $i=1,2,\ldots$, the $i$-th outputted
bit depends only on bits $1,2,\ldots,i$ of the input string.  According to
the results of this paper (see Proposition \ref{Bool}), 
the transformation of $\mathbb Z_2$ performed by the automaton is compatible,
that is, satisfy the 2-adic Lipschitz condition with a constant
1 and thus is continuous.  So 2-adic analysis can probably be of use in automata
theory.

\section*{Acknowledgement} I thank Branko Dragovich, Franco Vivaldi, and Igor Volovich 
for 
their interest to my work.
My special thanks to Andrei Khrennikov 
for encouraging discussions and hospitality during my stay at  V\"axj\"o university. 

%
%


\end{document}